\documentclass[twoside,12pt]{article}

\usepackage{amsmath}
\usepackage{amssymb}
\usepackage{amscd}
\usepackage{amsthm}

\numberwithin{equation}{section}

\theoremstyle{plain}

\theoremstyle{remark}

\theoremstyle{definition}

\newcommand{\Z}{\mathbb Z}
\newcommand{\R}{\mathbb R}
\newcommand{\C}{\mathbb C}
\newcommand{\h}{\mathcal H}
\newcommand{\Q}{\mathbb Q}

\newcommand{\n}{\mathcal N}

\def\ra{\rightarrow}

\def\e{\emph}
\def\i{\infty}

\def\b{\begin}
\newcommand{\ol}{\overline}

\begin{document}

\title{
{Quasiconformal maps on non-rigid Carnot groups}}
\author{Xiangdong Xie\footnote{Partially supported by NSF grant DMS--1265735.}}
\date{  }

\maketitle

\begin{abstract}
 We study quasiconformal maps on non-rigid   Carnot groups  equipped with Carnot metric.
We show that   for most non-rigid Carnot groups $N$,   all
quasiconformal maps on  $N$   must be biLipschitz.

\end{abstract}




 {\small {\bf{Mathematics Subject
Classification (2000).}}  22E25,   30L10,   53C17.   







\setcounter{section}{0}
  \setcounter{subsection}{0}

\section{Introduction}\label{s0}

We study quasiconformal maps on non-rigid   Carnot groups  equipped with Carnot metric.    We show that even most
non-rigid Carnot groups are rigid with respect to quasiconformal map. To be more precise,
  for most non-rigid Carnot groups  $N$,  all
quasiconformal maps on  $N$   must be biLipschitz.

A Carnot group  $N$ is called \e{$C^2$-rigid}  if the space of   $C^2$ contact maps
  between domains in   $N$ is finite dimensional,  and is called \e{non-rigid} otherwise.   In particular,  the space of $C^2$ quasiconformal maps on rigid Carnot groups  equipped with Carnot  metric   is finite dimensional.    On the other hand,  if   $\n=V_1\oplus \cdots \oplus
  V_r$ is the Lie algebra of $N$  and there exists some   nonzero   
    $X\in V_1$ such that the linear map
 $ad\; X:  \n\ra\n$, $ad\; X(Y)=[X,Y]$ has rank at most one, then  
    the space of biLipschitz maps from $N$ to $N$ is infinite dimensional (see   Proposition \ref{examplebilip}).
   Ottazzi  and      Warhurst (\cite{OW}, Theorem 1)  showed that a  Carnot group is non-rigid if and only if the complexification  $\n\otimes \C$ of the Lie algebra $\n$   of $N$ has
 an element of rank at most one.

Let $K\ge 1$ and $C>0$. A bijection $F:X\ra Y$ between two
     metric spaces is called a $(K, C)$-\e{quasi-similarity}  if
\[
   \frac{C}{K}\, d(x,y)\le d(F(x), F(y))\le C\,K\, d(x,y)
\]
for all $x,y \in X$.

Clearly a map is a quasi-similarity if and only if it is biLipschitz.
 The point here is that often there is control on $K$  but not on $C$.  In this case, the notion of quasi-similarity provides more information about the distortion.

We say that a Carnot group  $N$   is {\it quasisymmetrically rigid} if
 every $\eta$-quasisymmetric map from $N$ to $N$   is a $(K, C)$-quasi-similarity,  where
 $N$ is equipped with a  Carnot metric and  $K$ is a constant depending only on $\eta$.   See Section \ref{maps} for the definition of  quasisymmetric map.
  A  Carnot algebra is called quasisymmetrically rigid if the corresponding Carnot group is quasisymmetrically rigid.
  The main result of this paper says that most non-rigid  Carnot groups are
quasisymmetrically rigid.     We next describe those exceptional non-rigid Carnot groups.

 The first class of exceptional non-rigid Carnot groups   are 
  suitable quotients of direct products of the same Heisenberg group.

Let   $m, n\ge 1$ be integers.     Let $\tilde \n$ be the direct sum of $n$ copies of  the 
  $m$-th  Heisenberg algebra
 ${\mathcal H}_\R^m$.
 In other words,  $\tilde \n=\tilde {\mathcal H}_1\oplus \cdots \oplus \tilde
 {\mathcal H}_n$, where each $\tilde {\mathcal H}_j={\mathcal H}_\R^m$.
  Let $\tilde V_{1,j}$ and  $\tilde L_j$ respectively  be the first and second layers of   $\tilde {\mathcal H}_j$.  Let $\tilde V_1$ and $\tilde V_2$ respectively  be the first and second layers of  $\tilde \n$.
  Then     $\tilde V_1=\tilde V_{1,1}\oplus \cdots \oplus \tilde V_{1,n}$
  and $\tilde V_2=\tilde L_1\oplus \cdots \oplus \tilde L_n$.
  Since $\tilde V_2$ is central in  $\tilde \n$,  every  linear subspace $V\subset \tilde V_2$ is an ideal of $\tilde \n$ and  $\tilde \n/V=\tilde V_1\oplus(\tilde V_2/V)$ is  a   Carnot algebra.
Let  $V\subset  \tilde V_2$ be a linear subspace and $G_2\subset GL(\tilde V_2)$ be a
  group of linear transformations of $\tilde V_2$.
  We call $\n:=\tilde \n/V$  a  \e{Heisenberg product algebra}  if the following conditions
 are satisfied:\newline
  (1)  $\tilde L_j\cap V=\{0\}$ for all $j$;\newline
  (2)   $(V+\tilde L_i)\cap (V+\tilde L_j)=V$ for all $i\not=j$;\newline
  (3)  $V$ is $G_2$-invariant; that is,  $g_2(V)=V$ for all $g_2\in G_2$;\newline
  (4)   $G_2$ permutes the set  $\{\tilde L_1, \cdots,  \tilde L_n\}$ and acts transitively on the set.\newline
 A  Carnot group is called a  \e{Heisenberg product   group}  if its Lie algebra is a
Heisenberg product algebra.

 The  second  class of exceptional non-rigid Carnot groups   are 
  suitable quotients of direct products of the same   complex  Heisenberg group.

Recall that,  the 
$m$-th   ($m\ge 1$)  complex  Heisenberg algebra
 ${\mathcal H}_\C^m=\C^{2m}\oplus \C$ is a   complex Lie algebra of dimension $2m+1$.
  If   $e_1, \cdots, e_{2m}, \eta$ are the standard basis of $\C^{2m}\oplus \C$,
 then the only non-trivial bracket relations (over $\C$)  are $[e_i, \, e_{m+i}]=\eta$,  $1\le i\le m$.

Let   $m, n\ge 1$ be integers.     Let $ \tilde \n$ be the direct sum of $n$ copies of
the 
  $m$-th  complex  Heisenberg algebra
 ${\mathcal H}_\C^m$.
 In other words,  $\tilde \n=\tilde {\mathcal H}_1\oplus \cdots \oplus \tilde
 {\mathcal H}_n$, where each $\tilde {\mathcal H}_j={\mathcal H}_\C^m=\C^{2m}\oplus \C$. 
  Let $\tilde V_{1,j}$ and  $\tilde L_j$ respectively  be the first and second layers of   $\tilde {\mathcal H}_j$. 
We fix a graded isomorphism
 $f_j:  {\mathcal H}_\C^m=\C^{2m}\oplus \C\ra  \tilde {\mathcal H}_j$.   Then $f_j|_{\{0\}\oplus\C}:   \C\ra \tilde L_j$
  is a real linear isomorphism.  
 Let $\tilde V_1$ and $\tilde V_2$ respectively  be the first and second layers of  $\tilde \n$.
  Then     $\tilde V_1=\tilde V_{1,1}\oplus \cdots \oplus \tilde V_{1,n}$
  and $\tilde V_2=\tilde L_1\oplus \cdots \oplus \tilde L_n$.
  Since $\tilde V_2$ is central in  $\tilde \n$,  every  real  linear subspace $V\subset \tilde V_2$ is an ideal of $\tilde \n$ and  $\tilde \n/V=\tilde V_1\oplus(\tilde V_2/V)$ is  a   Carnot algebra.
Let  $V\subset  \tilde V_2$ be a  real  linear subspace and $G_2\subset GL(\tilde V_2)$ be a
  group of   real  linear transformations of $\tilde V_2$.
  We call $\n:=\tilde \n/V$  a  \e{complex  Heisenberg product algebra}  if the following conditions
 are satisfied:\newline
  (1)  $\tilde L_j\cap V=\{0\}$ for all $j$;\newline
  (2)  $V$ is $G_2$-invariant; that is,  $g_2(V)=V$ for all $g_2\in G_2$;\newline
  (3)   $G_2$ permutes the set  $\{\tilde L_1, \cdots,  \tilde L_n\}$ and acts transitively on the set;\newline
 (4) For each $g_2\in G_2$ and every $j$, if $g_2(\tilde L_j)=\tilde L_{\sigma(j)}$
  for some    $1\le \sigma(j)\le n$, then  there is some $0\not=a_j\in\C$ such that   the map  $g_{2,j}:=f^{-1}_{\sigma(j)}\circ g_2\circ f_j:\C\ra \C$   has exactly one of the two forms
  $g_{2,j}(z)=a_jz$ or  $g_{2,j}(z)=a_j\bar z$,  where $\bar z$ is the   
 complex conjugate   of $z$.  
   In other words,  $g_2|_{\tilde L_j}$ is a similarity.      
  \newline
A  Carnot group is called a  \e{complex  Heisenberg product   group}  if its Lie algebra is a
complex  Heisenberg product algebra.

Here is the main result of the paper:

\b{Th}\label{main}
Let $N$ be a non-rigid  Carnot group. If $N$ is not one of the following  three classes of groups, then
 it is quasisymmetrically rigid:\newline
 (1)      Euclidean groups;\newline
 (2)  Heisenberg product   groups;\newline
  (3)  complex  Heisenberg product   groups.
\end{Th}

It is known that there exist non-biLipschitz quasisymmetric  maps
 on Euclidean groups  \cite{GV}  and Heisenberg groups \cite{B}.
  It is an open question  whether the other Heisenberg product   groups and 
complex  Heisenberg product   groups are quasisymmetrically rigid. 
Two prominant examples among these groups are complex Heisenberg groups and the direct product of at least two copies of the same Heisenberg group.


The results in this paper extend previous results by the author
 for reducible Carnot groups  \cite{X3},   model Filiform groups   \cite{X1}  and
 $2$-step Carnot groups with reducible first layer \cite{X2}.

The first rigidity theorem about quasiconformal maps of Carnot groups is due to  Pansu.  He proved  that (\cite{P})   every quasiconformal map
  of the quarternionic Heisenberg group is a composition of left translations and graded automorphisms.
   In particular, the space of quasiconformal maps
 is finite dimensional.  
 Our results are of somehow different nature.  For example,   Proposition \ref{examplebilip}
    says  that the space of biLipschitz maps of Carnot groups  of rank  one
 is infinite dimensional.  Our main result  implies  that most of these  groups are
quasisymmetrically rigid.



  The results in this paper have implications for the large scale geometry of
 negatively curved
 homogeneous manifolds. 
  Each Carnot group arises  as  the (one point complement  of) ideal boundary  of some negatively
  curved   homogeneous manifold  \cite{H}.
   Furthermore,  each quasiisometry of the negatively curved homogeneous manifold associated to a  quasisymmetrically  rigid Carnot group  
   is a rough isometry, that is, it must preserve the distance up to an additive constant.

In Section \ref{prelimi}  we  collect definitions and results that shall be needed later.
  In Section \ref{invariant}  we show that  high step non-rigid Carnot algebras have reducible first layer.  In Section \ref{2step}  we give a   characterization of non-rigid Carnot algebras  with
irreducible first  layer.  
     In Section \ref{highstep}  we prove the main result (Theorem \ref{main}) of the paper.
Finally in Section \ref{examples} we construct an infinite dimensional space of
biLipschitz maps on  rank one Carnot groups.

\noindent {\bf{Acknowledgment}}. {This work was   initiated   while the author was attending the workshop
\lq\lq Interactions between analysis and geometry" at IPAM,  University of California at Los Angeles  from March to June 2013.  I would like to thank IPAM  for   financial support, excellent working conditions and
   conducive atmosphere.
 I also would like to thank David Freeman and Tullia Dymarz for discussions about Carnot groups.  }

\section{Preliminaries}\label{prelimi}


In this Section we collect definitions and results that shall be needed later.

\subsection{Carnot groups  and Carnot algebras}\label{basics}

A   \e{Carnot Lie algebra} is a finite dimensional Lie algebra
$\mathcal G$   together with  a direct sum   decomposition
    $\mathcal
G=V_1\oplus V_2\oplus\cdots \oplus V_r$
  of   non-trivial   vector subspaces
 such that $[V_1,
V_i]=V_{i+1}$ for all $1\le i\le r$,
    where we set $V_{r+1}=\{0\}$.  The integer $r$ is called the
    degree of nilpotency of $\mathcal
G$. Every Carnot algebra
 $\mathcal
G=V_1\oplus V_2\oplus\cdots \oplus V_r$  admits a one-parameter
 family of automorphisms $\lambda_t: \mathcal
G \ra \mathcal G$, $t\in (0, \i)$, where
 $\lambda_t(x)=t^i x$ for  $x\in V_i$.
  Let   $\mathcal
G=V_1\oplus V_2\oplus\cdots \oplus V_r$
    and $\mathcal
G'=V'_1\oplus V'_2\oplus\cdots \oplus V'_s$  be two  Carnot
    algebras.
  A Lie algebra  homomorphism
  $\phi: \mathcal
G\ra \mathcal G'$
     is graded if $\phi$ commutes with $\lambda_t$ for
  all $t>0$; that is, if $\phi\circ \lambda_t=\lambda_t\circ
  \phi$.  We observe that $\phi(V_i)\subset V'_i$ for all $1\le i\le
  r$.

A  simply connected nilpotent Lie group is a  \e{Carnot group}
 if its Lie algebra is a Carnot algebra.
    Let $G$ be a Carnot group with Lie algebra
      $\mathcal G=V_1\oplus \cdots \oplus V_r$.  The subspace $V_1$ defines
      a left invariant distribution  $H G\subset TG$ on $G$.    We fix a left invariant inner product on
          $HG$.
           An
      absolutely continuous curve $\gamma$ in $G$  whose velocity vector
       $\gamma'(t)$  is contained in  $H_{\gamma(t)} G$ for a.e. $t$
        is called  a horizontal curve.
          By Chow's theorem ([BR, Theorem 2.4]),   any two points
  of $G$ can be  connected by horizontal curves. Let $p,q\in G$, the
  \e{Carnot   metric} $d_c(p,q)$   between them is defined as
  the infimum of length of horizontal curves that join $p$ and $q$.

  Since the inner product on   $HG$ is left invariant, the Carnot
  metric on $G$ is also left invariant.  Different choices of inner
  product on $HG$ result in Carnot metrics that are biLipschitz
  equivalent.
    The Hausdorff dimension of $G$ with respect to  a  Carnot metric
    is  given by $\sum_{i=1}^r i\cdot \dim(V_i)$.

Recall that, for a simply connected nilpotent Lie group $G$ with Lie
algebra $\mathcal G$, the exponential map
  $\exp: {\mathcal G}\ra G$ is a diffeomorphism.  Under this identification the Lesbegue   measure on
  $\mathcal G$  is a  Haar measure on $G$.
 Furthermore, the
  exponential map induces 
  a  one-to-one correspondence between
    Lie subalgebras of $\mathcal G$   and
  connected Lie subgroups of $G$.

 It is often   more  convenient to work with homogeneous distances defined using norms
 than with Carnot metrics.  Let
$\mathcal
G=V_1\oplus V_2\oplus\cdots \oplus V_r$  be a  Carnot algebra.
  Write $x\in \mathcal G$ as $x=x_1+\cdots+ x_r$ with $x_i\in V_i$.
  Fix a norm   $|\cdot|$ on each layer.  Define a norm $||\cdot||$ on $\mathcal G$
   by:
$$||x||=\sum_{i=1}^r |x_i|^{\frac{1}{i}}.$$
 Now define a homogeneous distance on $G=\mathcal G$  by:  $d(g,h)=||(-g)*h||$.
  In general,   $d$ is only a quasimetric.  However,  
$d$ and $d_c$ are always  biLipschitz equivalent.   That is, there is a constant $C\ge 1$ such that    $d(p,q)/C\le d_c(p,q)\le C\cdot   d(p,q)$ for all $p, q\in G$.
  It is often possible to  calculate  or estimate  $d$  by using  the BCH formula (see  Subsection \ref{BCH}).
     Since we are only concerned with quasiconformal maps and biLipschitz maps,
       it does not matter whether we use $d$ or $d_c$.


\subsection{The  Baker-Campbell-Hausdorff  formula}\label{BCH}

Let  $G$ be a simply connected nilpotent Lie group with Lie algebra $\mathcal{G}$.
  The exponential map  $\text{exp}:  {\mathcal{G}}\ra G$ is a diffeomorphism.
   One can then
     pull back the group operation from $G$ to get a group stucture
 on $\mathcal{G}$.     This group structure can be described by the   Baker-Campbell-Hausdorff formula
 (BCH formula in short),  which expresses the product $X*Y$ ($X, Y\in {\mathcal{G}}$)
    in terms of the iterated Lie brackets  of $X$ and $Y$.
     The group operation in $G$ will be denoted by $\cdot$.
   The pull-back group operation   $*$  on $\mathcal{G}$ is defined as follows.
      For $X,   Y\in \mathcal{G}$,   define
   $$X*Y=\text{exp}^{-1}(\text{exp} X\cdot \text{exp} Y).$$
  Then the BCH formula (\cite{CG},    page 11)  says
  $$X*Y=\sum_{n>0}\frac{(-1)^{n+1}}{n}\sum_{p_i+q_i>0, 1\le i\le n}\frac{(\sum^n_{i=1}(p_i+q_i))^{-1}}
  {p_1!q_1!\cdots p_n!q_n!} (\text{ad} X)^{p_1}(\text{ad} Y)^{q_1}\cdots   (\text{ad} X)^{p_n}(\text{ad} Y)^{q_n-1}Y$$
  where $\text{ad}\, A (B)=[A,   B]$.     If  $q_n=0$,  the term in the sum is  $\cdots (\text{ad} X)^{p_n-1}X$;
  if  $q_n>1$  or if $q_n=0$  and $p_n>1$, then the term  is zero.
  The first a few terms are well known,
\b{align*}
 X*Y& =X+Y+\frac{1}{2}[X,Y]+\frac{1}{12}[X,[X,Y]]
-\frac{1}{12}[Y, [X, Y]]\\
& -\frac{1}{48}[Y,[X,[X,Y]]]-\frac{1}{48}[X,[Y,[X,Y]]]+(\text{commutators  in five or more terms}).
\end{align*}

\subsection{Quasisymmetric    map}\label{maps}

Here we recall the definitions of   quasisymmetric maps  and pointwise Lipschitz constant.



Let $\eta: [0,\i)\ra [0,\i)$ be a homeomorphism.
    A homeomorphism
$F:X\to Y$ between two metric spaces is
\e{$\eta$-quasisymmetric} if for all distinct triples $x,y,z\in X$,
we have
\[
   \frac{d(F(x), F(y))}{d(F(x), F(z))}\le \eta\left(\frac{d(x,y)}{d(x,z)}\right).
\]
  If $F: X\ra Y$ is an $\eta$-quasisymmetry, then
  $F^{-1}:    Y\ra X$ is an $\eta_1$-quasisymmetry, where
$\eta_1(t)=(\eta^{-1}(t^{-1}))^{-1}$. See \cite{V}, Theorem 6.3.
 A    homeomorphism between metric spaces
  is quasisymmetric if it is $\eta$-quasisymmetric for some $\eta$.

We remark that
     quasisymmetric homeomorphisms between general metric spaces
  are quasiconformal. In the case of Carnot groups  (and more
  generally  Loewner  spaces),  a   homeomorphism  is quasisymmetric if and only
  if it is quasiconformal, see \cite{HK}.



Let $F: X\ra  Y$ be a   homeomorphism   between two
metric spaces.
We define for every $x\in X$  and $r>0$,
\begin{align*}
   L_F(x,r)&=\sup\{d(F(x), F(x')):   d(x,x')\le r\},\\
   l_F(x,r)&=\inf\{d(F(x), F(x')):   d(x,x')\ge r\},
\end{align*}
and set
\[
   L_F(x)=\limsup_{r\ra 0}\frac{L_F(x,r)}{r}, \ \
   l_F(x)=\liminf_{r\ra 0}\frac{l_F(x,r)}{r}.
\]
  Then
\b{equation}\label{e3.001}
  L_{F^{-1}}(F(x))=\frac{1}{l_F(x)} \ \text{ and }\ l_{F^{-1}}(F(x))=\frac{1}{L_F(x)}
\end{equation}
for any $x\in X$. If $F$ is an $\eta$-quasisymmetry, then
\b{equation}\label{e3.0}
L_F(x,r)\le \eta(1)l_F(x, r)
\end{equation}
 for all $x\in X$ and $r>0$. 

\subsection{Previous results}\label{previous}

Here we collect some results that will be used  in this paper. 

A   $C^2$ diffeomorphism between domains of   Carnot groups  whose differential preserves the horizontal bundle is called a contact map.   A  Carnot group   $N$ is  
\e{rigid} if the space of   $C^2$-contact maps between domains of $N$ is finite
    dimensional, otherwise,  $N$ is    \e{non-rigid}.   A Carnot algebra is called rigid  if the corresponding Carnot
 group is rigid. Similarly  a  Carnot algebra is called non-rigid  if the corresponding Carnot
 group is non-rigid. 

Non-rigid Carnot groups (algebras) can be  characterized in 
terms of rank.  
  The rank of an element in a Lie algebra     appeared in \cite{O}.

For an element $x\in \mathcal{N}$ in a  Lie algebra, let $\text{rank}(x)$ be the rank of
 the linear transformation  $\text{ad}(x):  \mathcal{N}\ra   \mathcal{N}$,    $\text{ad}(x)(y)=[x,y]$.   In other words,
  $\text{rank}(x)$  is the dimension of the image of $\text{ad}(x)$.
 It is easy to see that $\text{rank}(x)=\text{rank}(F(x))$  for any isomorphism  $F: \n_1\ra \n_2$.


\b{Th}\label{nonrch}\e{(\cite{OW},  Theorem 1)}
  Let $\n$  be a  Carnot algebra.  Then  $\n$ is non-rigid if and only if
  there exists a   nonzero vector   $X$  in the first layer of the complexification of $\n$   such that  $\text{rank}(X)\le 1$.

\end{Th}

We say that a Carnot group  $N$   is {\it quasisymmetrically rigid} if  
 every $\eta$-quasisymmetric map from $N$ to $N$  is a $(K, C)$-quasi-similarity,   where 
   $N$ is  equipped with a  Carnot metric and  $K$ is a constant depending only on $\eta$.  
 A Carnot algebra is 
quasisymmetrically rigid  if the corresponding Carnot group is.

Let $\text{Aut}_g(\mathcal N)$ be the group of graded isomorphisms of a  Carnot algebra $\mathcal N$. We say $V_1$ is 
{\it reducible} (or the first layer of $\mathcal N$ is reducible) if  there is a non-trivial proper
 linear subspace $W\subset V_1$ such that $A(W)=W$ for every 
 $A\in \text{Aut}_g(\mathcal N)$.

\b{Th}\label{2step1}  \e{(\cite{X2}, Theorem 1.1)}
  Let $\n$ be a $2$-step Carnot  algebra. If the first layer of  $\n$  is reducible,
  then $\n$ is quasisymmetrically rigid.
\end{Th}

\b{Th}\label{2step2} \e{(\cite{X2}, Theorem 1.2)}
Suppose $W_1\subset V_1$ is a non-trivial proper subspace that is invariant under the action of
 $\text{Aut}_g(\mathcal{N})$. 
  If   there is some $X\in V_1\backslash W_1$ such that $[X, W_1]\subset [W_1, W_1]$,
  then $\n$ is quasisymmetrically rigid.
\end{Th}

 The following result is very useful for induction argument:

\b{Th}\label{2step3} \e{(\cite{X2}, Theorem 1.3)}
 Suppose $W_1\subset V_1$ is a non-trivial proper subspace that is invariant under the action of
 $\text{Aut}_g(\mathcal{N})$.  
    If  the Lie subalgebra  $<W_1>$   generated by $W_1$ is quasisymmetrically rigid, then so is
$\n$.  
\end{Th}

We shall also need the following result (Section \ref{abelian}).

\b{Prop}\label{p3} \e{(\cite{X3}, Proposition   3.4) }
 Let  $G$  and $G'$ be two Carnot groups
  with Lie algebras  $\mathcal G=V_1\oplus \cdots \oplus V_m$  and $\mathcal G'=V'_1\oplus \cdots \oplus  V'_n$ respectively.    Let     $W_1\subset V_1$,  $W'_1\subset V'_1$ be  linear
  subspaces. Denote by $\mathcal W\subset \mathcal{G}$ and
$\mathcal {W'}\subset \mathcal{G}'$   respectively the Lie subalgebras
  generated by $W_1$ and $W'_1$.    Let $W\subset G$ and $W'\subset G'$ respectively  be the
  connected   Lie
  subgroups of $G$ and $G'$    corresponding to
  $\mathcal W$  and    $\mathcal {W'}$.
 Let $F: G\ra G'$ be a quasisymmetric  homeomorphism.
  If     $dF(x)(W_1)\subset W'_1$ for a.e. $x\in G$,
     then
   $F$ sends left cosets of $W$ into left  cosets of $W'$.

  \end{Prop}

\section{High step non-rigid Carnot algebras have reducible first layer}\label{invariant}

The goal of this section is to show that if $\n$ is non-rigid and has step $r\ge 3$, then
    the first layer is reducible.


  Let $\mathcal N=V_1\oplus \cdots \oplus V_r$ be a Carnot algebra  over a field $K$. 
 We are only interested in the cases $K=\R$ and $K=\C$.
   Define
  $$r_{1,K}(\mathcal N)=\min\{\text{rank}(x):  0\not=x\in  V_1\}. $$
  When $K=\R$,  we simply write $r_1(\n)$ for $r_{1, \R}(\n)$ and call it 
the {\it rank} of $\mathcal N$ and of the 
    corresponding Carnot group $N$. 
  Let $W_1\subset  V_1$ be the linear subspace of $V_1$ spanned by elements 
 $x\in V_1$ satisfying $\text{rank}(x)=r_1(\mathcal N)$.  
  Since graded isomorphisms preserve  the first layer $V_1$,  we see that $W_1$ is a 
 non-trivial subspace of $V_1$ invariant under the action of 
 $\text{Aut}_g(\mathcal N)$.     
   The question  is when $W_1$ is a proper subspace of $V_1$.

\b{Le}\label{rank0}
Let $\n$ be a non-abelian Carnot algebra.   If  $r_1(\n)=0$, then $W_1$ is a proper subspace of $V_1$.
\end{Le}

\b{proof}
Let $X, Y\in V_1$ with $\text{rank}(X)=\text{rank}(Y)=0$, and $a, b\in \R$.  
 Then $[X, \n]=[Y, \n]=0$. Hence $[aX+bY, \n]=a[X, \n]+b[Y, \n]=0$.  It follows that $W_1$ is the set of elements in $V_1$ with rank $0$. Since $\n$ is non-abelian, $W_1\not=V_1$.

\end{proof}

\b{Le}\label{rank0.1}
  For any Carnot  Lie algebra $\mathcal N$,  $r_1(\n)=0$
 if and only if
   $\mathcal{N}$ can be written as a direct sum of an Euclidean algebra  and another
 Carnot algebra.

\end{Le}

\b{proof}
    First  suppose $\mathcal{N}=\R^k\oplus \mathcal{N}'=(\R^k\oplus V'_1)\oplus V'_2\oplus  \cdots  \oplus V'_r$ is a
 direct sum of an Euclidean algebra  and another
 Carnot algebra.
Clearly  every $X\in \R^k\backslash\{0\}\subset V_1=\R^k\oplus V'_1$ has rank $0$.

  Conversely,  assume
   there exists a nonzero $X\in V_1$ with $\text{rank}(X)=0$.
    Let $V'_1\subset V_1$ be a  codimension $1$ subspace complementary to $\R X$.
  Set $\mathcal{N}'=V'_1\oplus V_2\oplus \cdots \oplus V_r$.  Then it is easy to check that $[V'_1, V'_1]=V_2$ and $[V'_1, V_i]=V_{i+1}$ for $2\le i\le r$, where $V_{r+1}=0$.
   It follows that  $\mathcal{N}'$  is also a Carnot algebra.     Now it is easy to check that the map  
 $\R\oplus   \mathcal{N}'  \ra  \mathcal{N}$,  $(a,Y)\ra   aX+Y$,  is a
graded  isomorphism of   Carnot algebras.

\end{proof}

\b{Le}\label{rank1}
Suppose  $r_1({\mathcal{N}})=1$  and 
    $W_1=V_1$.
  Then   $\mathcal{N}$ is a $2$-step Carnot algebra.

\end{Le}

\b{proof}
The assumptions imply that    there is a vector space basis $X_1, \cdots, X_m$ of $V_1$ satisfying
 $\text{rank}(X_i)=1$ for all $1\le i\le m$.  Since $\n$  is generated by $V_1$, we must have
 $[X_i, V_1]\not=0$. Hence   for each $i$  there is some $j_i$ such that $Y_i:=[X_i, X_{j_i}]\not=0$.
      Notice that $Y_i\in V_2$.  Hence
   $V_2\not=0$.
Since $\text{rank}(X_i)=1$,
  we have $[X_i,  \n]=\R Y_i$.   
  Since $[X_i,  V_j]\subset V_{j+1}$,
    we have   $[X_i, V_j]=0$ for  all $j\ge 2$.
In particular,  $[X_i, V_2]=0$.    Since
$X_1, \cdots, X_m$   form a basis of $V_1$, we have $V_3=[V_1, V_2]=0$.  So $\mathcal{N}$ is a $2$-step Carnot algebra.

\end{proof}

\b{Cor}\label{highrank}
  Let $\mathcal{N}$ be a $r$-step Carnot algebra with $r\ge 3$.
  If  $r_1(\n)\le 1$,
then
$W_1\not=V_1$.

\end{Cor}

\b{proof}
If $r_1({\mathcal{N}})=0$,  the claim follows from Lemma \ref{rank0}. 
 So we  assume
   $r_1({\mathcal{N}})=1$.    Since $\mathcal{N}$ is $r$-step with $r\ge 3$, the
 Corollary now follows from   Lemma \ref{rank1}.

\end{proof}

Let $\mathcal{N}_\C=\mathcal{N}\otimes \C=(V_1\otimes \C)\oplus \cdots \oplus (V_r\otimes \C)$ be
  the complexification of  $\mathcal{N}$.
Let  $r_{1, \C}(\n) =\min\{\text{rank}(X)|  0\not=X\in V_1\otimes \C\}$  and
   $W_{1, \C}\subset V_1\otimes \C$ be the   complex linear subspace  of
  $V_1\otimes \C$
spanned by elements $0\not=X\in  V_1\otimes \C$ with
  $\text{rank}(X)=r_{1, \C}(\n)$.    

\b{Le}\label{crank1}
Suppose $r_1({\mathcal{N}})\ge 2$.  Then
$r_{1, \C}(\n) \ge 1$.    If   $r_{1, \C}(\n)= 1$  and
$W_{1, \C}=V_1\otimes \C$,  then $\mathcal{N}$ is a  $2$-step
  Carnot algebra.

\end{Le}

\b{proof}
Suppose there is some $0\not=X=X_1+i X_2\in V_1\otimes \C$ with
$\text{rank}(X)=0$, where $X_1, X_2\in V_1$.
  Then $[X, \n]\subset [X,  \n_\C]=0$.  
  It follows that 
$[X_1, \n]=0=[X_2, \n]$.
 This implies   there is a nonzero element  ($X_1$  or $X_2$)  in $V_1$
 with rank $0$, contradicting the assumption
   $r_1({\mathcal{N}})\ge 2$.  Hence $r_{1,  \C}(\n)\ge 1$.

 Now suppose
   $r_{1,  \C}(\n)=1$   and
$W_{1, \C}=V_1\otimes \C$.  Then the proof of Lemma \ref{rank1}   shows
  $[V_1\otimes \C, V_2\otimes \C]=0$.
  Hence  $[V_1\otimes \C, V_2\otimes \C]\supset [V_1, V_2]=0$, which implies $\mathcal{N}$ is $2$-step.

\end{proof}

\b{Le}\label{l4.12}
If  $r_{1, \C}(\n)\le 1$, then $r_1(\n)\le 2$.
\end{Le}

\b{proof}
In the proof of Lemma \ref{crank1} we already showed that if
 $r_{1, \C}(\n)=0$, then  $r_1(\n)=0$.    So we shall assume
$r_{1, \C}(\n)= 1$.     Let $0\not=X=X_1+iX_2\in V_1\otimes \C$ $(X_1, X_2\in V_1$)
  with   $\text{rank}(X)=1$.
  By multiplying $i$ if necessary we may assume $X_1\not=0$.
    We must have $[X, V_1]\not=0$, otherwise  $[X, \n \otimes \C]=0$.
  So there is some $Y_0\in V_1$ such that     $[X, Y_0]\not=0$.  Since
$r_{1, \C}(\n)= 1$, we have $[X, \n \otimes \C]=\C [X, Y_0]$.
  In particular, for any $Y\in \n$,  there are $a, b\in \R$ such that
 $$[X_1+iX_2,  Y]=[X, Y]=(a+ib)[X, Y_0]=(a+ib)[X_1+iX_2,  Y_0].$$
Comparing   the real parts of both sides, we obtain  $[X_1, Y]=a[X_1, Y_0]-b[X_2, Y_0]$.
 This implies  $[X_1, \n]\subset \R[X_1, Y_0]+\R [X_2, Y_0]$.
 Hence $\text{rank}(X_1)\le 2$.

\end{proof}

\b{Cor}\label{crank1cor}
 Let   $\hat{W_1}=\{X\in V_1|  \text{there exists}\; Y\in V_1 \;\text{such that }\;\;   X+iY\in W_{1,\C}\}$.
If     $r_1({\mathcal{N}})\ge 2$,  $r_{1, \C}(\n)  = 1$  and  $\mathcal{N}$ is $r$-step with $r \ge 3$,  then
   $\hat{W_1}\not=V_1$.

\end{Cor}

\b{proof}
  By Lemma \ref{crank1}  and our assumption that $r\ge 3$,      we have $W_{1, \C}\not=V_1\otimes \C$.  Let $X_1, \cdots, X_m$ be a basis for the complex vector space
   $W_{1, \C}$  satisfying  $\text{rank}(X_j)=1$.
    We write $X_j=Y_j+i   Z_j$ with $Y_j,  Z_j\in V_1$.  Then it is   easy  to check that
 as a real vector space  $\hat{W_1}$ is spanned by $Y_1, Z_1, \cdots, Y_m, Z_m$.
Suppose $\hat{W_1}=V_1$.  Then
$Y_1, Z_1, \cdots, Y_m, Z_m$  span  $V_1$.    As seen in the proof of
  Lemma \ref{rank1},  we have $[X_j, V_2\otimes \C]=0$.  In particular,
   $[Y_j+i Z_j, V_2]=[X_j, V_2]=0$.   It follows that $[Y_j, V_2]=0=[Z_j, V_2]$.
  Since
$Y_1, Z_1, \cdots, Y_m, Z_m$  span  $V_1$,   we have   $V_3=[V_1, V_2]=0$,
 contradicting our assumption that $r\ge 3$.
  Hence  $\hat{W_1}\not=V_1$.

\end{proof}

\b{Prop}\label{invariantspace}
  If  a Carnot algebra $\mathcal{N}$ has step $r\ge 3$ and is non-rigid, then
 $V_1$  is reducible.

\end{Prop}

\b{proof}  First of all,   $\text{rank}(L(X))=\text{rank}(X)$ for any Lie algebra isomorphism
 $L: \mathcal{N}_1 \ra\mathcal{N}_2$ and any $X\in \mathcal{N}_1$.
  A graded isomorphism  between Carnot algebras  also preserves the first layer.  Hence   for any $L\in \text{Aut}_g(\mathcal{N})$
 and any $X\in V_1$, we have $L(X)\in V_1$  and
    $\text{rank}(L(X))=\text{rank}(X)$.   It follows that $L(W_1)=W_1$.
  Similarly,    $L(W_{1, \C})=W_{1, \C}$  and $L(\hat {W_1})=\hat{W}_1$.
Hence $W_1$ and $\hat W_1$ are non-trivial linear subspaces of $V_1$ that are invariant under the action of  $\text{Aut}_g(\mathcal{N})$.   Next  we shall  decide whether they are proper subspaces of $V_1$.

 Since $\n$ is non-rigid,   Theorem \ref{nonrch}  implies
     $r_{1, \C}(\n)  \le 1$.
 If $r_1(\mathcal{N})\le 1$,  Corollary \ref{highrank}
  implies  $W_1$ is proper.
  Now assume $r_1(\mathcal{N})\ge 2$.
 Then Lemma  \ref{crank1}  and  Corollary \ref{crank1cor}  imply
  $\hat{W_1}\not=V_1$.

\end{proof}

\section{Non-rigid Carnot algebras with irreducible first  layer}\label{2step}

In this Section we give a   characterization of non-rigid Carnot algebras  with
irreducible first  layer.

Let   $\n=V_1\oplus \cdots \oplus V_r$  be a Carnot algebra over a field $K$.   We are only 
  interested in  the cases $K=\R$   and  $K=\C$. 
 Let $A_1=\{v\in V_1:  \text{rank}(v)=1\}$  be the set of elements in $V_1$ with rank $1$.
  Notice that if $v\in A_1$, then $[v,  \n]=K  Z$ for some $0\not= Z\in V_2$.
We define an equivalence relation on   $A_1$ as follows:
  $v\sim v'$ if    $[v, \n]=[v', \n]$.    It is easy to check that this is indeed an equivalence relation.

\b{Le}\label{rank1subspace}
Suppose $r_{1, K}(\mathcal N)=1$.
 Let  $A$ be an equivalence class. Then $A\cup \{0\}$  is a linear subspace of $V_1$.

\end{Le}

\b{proof}
$A\cup \{0\}$   is clearly closed under scalar multiplication.
   We need to show  that  it is also closed under addition. Let $X_1, X_2\in A\cup \{0\}$.
  We may assume $X_1,  X_2,  X_1+ X_2\not=0$.
  Then $X_1\sim X_2$ and so there is some $0\not= Z\in V_2$
 such that 
$[X_1,\, \n]=K Z=[X_2,\,  \n]$.  
So  $[X_1+ X_2,  \,\n]\subset K Z$  and   $\text{rank}(X_1+X_2)\le 1$.
  Since $X_1+X_2\not=0$  and $r_{1, K}(\mathcal N)=1$,  we must have
$\text{rank}(X_1+X_2)=1$  and
$[X_1+ X_2, \, \n]=K Z$.    Hence   $X_1+X_2\in A$.

\end{proof}

\b{Le}\label{l21}
  Let $A_1, A_2$ be two distinct equivalence classes. Then $[X_1, X_2]=0$ for any $X_1 \in  A_1$, $ X_2\in A_2$.

\end{Le}

\b{proof}
There are $0\not= Z_1, Z_2\in V_2$ such that $ K Z_1\not=K Z_2$ and
  $[X_1, \n]=K  Z_1$  and  $[X_2, \n]=K Z_2$.  So
 $[X_1, X_2]\in K Z_1\cap K Z_2=\{0\}$.

\end{proof}

We call a subset of the form $\tilde {A}:=A\cup \{0\}$ an extended equivalence class,
    where $A$ is an equivalence class.

\noindent
{\bf{Example}}\label{l1p}
Let $\mathcal N=V_1\oplus V_2$ be a $2$-step Carnot algebra defined as follows.
   The first layer
 $V_1$ has a vector space basis $X_1, X_2, Y$ and the second layer
  $V_2$ has a vector space basis $Z_1,  Z_2$.
  The only non-trivial  bracket  relations are $[X_1, Y]=Z_1$, $[X_2, Y]=Z_2$.
  It is easy to check that  $r_1(\n)=1$,   and that 
elements of rank $1$ in $V_1$ are exactly the nonzero elements in
 $\R X_1\oplus \R X_2$  and   $[a_1X_1+ a_2 X_2,\,   \n]=\R (a_1 Z_1+ a_2 Z_2)$.
   It follows that the extended equivalence classes are exactly the $1$-dimensional subspaces of
   $\R X_1\oplus \R X_2$.

 Consider  linear subspaces of $V_1$   that  consist  of  only   rank $1$ elements (and  $0$).   These subspaces are
partially ordered by inclusion.
 The above example shows that,
in general,  the extended equivalence classes are not maximal   with respect to inclusion    and
   their span is not a direct sum of  all the extended equivalence classes.
        This picture changes when   $V_1$ is spanned by rank $1$ elements.

\b{Le}\label{rankspan}
Let $K=\R$ or $\C$.    
    Let $\mathcal N$ be a Carnot algebra over $K$.
   Suppose $r_{1, K}(\mathcal N)=1$  and
  $V_1$  is spanned by rank $1$ elements.     Then\newline
 (1) $V_1$ is  the direct sum of extended equivalence classes;\newline
  (2)  Each extended equivalence class  is a maximal linear subspace (over $K$)  of $V_1$
  consisting of ($0$ and)  rank $1$ elements;\newline
  (3) For each   extended   equivalence  class $\tilde  A$,  the   Lie subalgebra
 $<\tilde  A>$  is isomorphic to a Heisenberg algebra over $K$.
\end{Le}

\b{proof}
   (1)  and (2).
 First we notice  that two distinct extended equivalence classes intersect trivially. Let
  $\tilde {A_1}$ and $\tilde {A_2}$ be two distinct extended equivalence classes.
 Since $A_1$ and $A_2$ are distinct equivalence classes, we have $A_1\cap A_2=\emptyset$.
  Hence   $\tilde {A_1}\cap \tilde {A_2}=\{0\}$.

Since every rank   $1$  element lies in an equivalence class, the assumption implies  that
   $V_1$ is spanned by extended equivalence classes.
    Hence  there are extended equivalence classes $\tilde {A_1}, \cdots, \tilde {A}_m$ such that
 $V_1=\tilde {A_1}+\cdots+ \tilde {A}_m$.   We may assume
  $$\tilde {A_1}+\cdots +\tilde {A_j}\not= \tilde {A_1}+\cdots +\tilde {A}_{j+1}$$
  for $1\le j\le m-1.$
  Now fix a subspace   $B_{j+1}\subset \tilde {A}_{j+1}$  such that   (1)
 $\tilde {A_1}+\cdots +\tilde {A}_{j+1}=\tilde {A_1}+\cdots +\tilde {A_j}+{B_{j+1}}$  and (2)
$(\tilde {A_1}+\cdots +\tilde {A_j})\cap B_{j+1}=\{0\}$.
  The first paragraph implies $B_2=\tilde {A_2}$.
  Set $B_1=\tilde {A_1}$.  
Then
  \b{equation}\label{e6.10}
V_1=B_1\oplus    B_2\oplus B_3\oplus \cdots \oplus  B_m.
\end{equation}
  Since $B_j\subset \tilde {A_j}$,  Lemma  \ref{l21}
 implies
 $[B_i, B_j]=0$ for $i\not=j$.
    Since $r_{1, K}(\mathcal N)=1$,
 for any $0\not=X\in B_j$, there exists some
$Y\in B_j$ with $[X, Y]\not=0$.

We claim that $\text{rank}(X)\ge 2$ for any  $X\in  V_1\backslash(\cup_j B_j)$.  Note that 
   (2) follows from the claim.  Now we prove the claim.
 Write $X=X_{j_1}+\cdots +X_{j_k}$  with $k\ge 2$,  $1\le j_1< \cdots <j_k\le m$,   $0\not=X_{j_i}\in B_{j_i}$.
    For each $1\le j\le m$,  let  $0\not=Z_j\in V_2$ be  such that
   $[Y,  \n]=K Z_j$ for every $Y\in A_j$.   By the preceding paragraph,
  there exists some $Y_{j_i}\in  B_{j_i}$  with $[X_{j_i}, Y_{j_i}]=a_{j_i} Z_{j_i}$ for some $a_{j_i}\not=0$.
  Now we have $[X, Y_{j_1}]=[X_{j_1}, Y_{j_1}]=a_{j_1} Z_{j_1}$  and $[X, Y_{j_2}]=[X_{j_2}, Y_{j_2}]=a_{j_2} Z_{j_2}$.
   Since $K Z_{j_1}\not=K Z_{j_2}$, we see that $\text{rank}(X)\ge 2$.
  This in  particular implies $\tilde{A}_3\cap (\tilde {A}_1\oplus \tilde {A}_2)=\{0\}$ because  the        elements in ${A_3}$ have rank
 $1$ and the elements in $(\tilde {A}_1\oplus \tilde{A}_2)\backslash (\tilde{A}_1\cup \tilde {A}_2)$ have rank  (at least) $2$.
 So we actually have $B_3=\tilde{A}_3$. By induction  we obtain $B_j=\tilde{A}_j$ for all $j$.
 Hence (1)   follows from (\ref{e6.10}).

(3)  For any $X\in   {A_j}$, there is some $Y\in \tilde {A}_j$ such that $[X,Y]=a Z_j$ for some $a\not=0$. Now (3) follows from the
  following Lemma.

\end{proof}

\b{Le}\label{heisenberg}
 Let $K$ be a field and
  $\mathcal  N=V_1\oplus V_2$ be a $2$-step Carnot algebra    over $K$  with  $\dim(V_2)=1$.  Suppose
   for any $0\not=X\in V_1$ there is some $Y\in V_1$ such that $[X, Y]\not=0$.
  Then  $\n$ is a  Heisenberg algebra over $K$.

\end{Le}

\b{proof}
  Let $Z\in V_2$ be a nonzero element.
Fix any $X_1\in V_1\backslash\{0\}$.   The assumption   implies  that  
    there is some $Y_1\in V_1$ with $[X_1, Y_1]=Z$.  Since  $\dim(V_2)=1$  and  $[V_1, V_1]\subset V_2$,  
   we see   that $\ker(\text{ad}\, X_1)\cap \ker(\text{ad} \,Y_1)$ has codimension $2$ in $V_1$  and
  $$V_1=K X_1\oplus K   Y_1\oplus (\ker(\text{ad} \,X_1)\cap \ker(\text{ad}\, Y_1)).$$
  If  $ \ker(\text{ad} \,X_1)\cap \ker(\text{ad}\, Y_1)=\{0\}$,  then   $\n$ is isomorphic to the first Heisenberg algebra over $K$  and we are done. 
Assume $ \ker(\text{ad} \,X_1)\cap \ker(\text{ad}\, Y_1)\not=\{0\}$.
  In this case, it is easy to check that the subalgebra $\n_1:=<(\ker(\text{ad} \,X_1)\cap 
\ker(\text{ad}\, Y_1))>$  satisfies the assumption in the Lemma.  
By induction we may assume 
$\n_1$ is the $m$-th Heisenberg algebra over $K$  for some $m\ge 1$.
  Hence there are $X_j, Y_j\in  \ker(\text{ad} \,X_1)\cap \ker(\text{ad}\, Y_1)$ for $2\le j\le m+1$       such that \newline
  (1)   $[X_j, Y_j]=Z$ for all $2\le j\le m+1$  and  $[X_i, X_j]=0$, $[X_i, Y_j]=0$, $[Y_i, Y_j]=0$ for $i\not=j$;\newline
 (2)  $\ker(\text{ad} \,X_1)\cap \ker(\text{ad}\, Y_1)=K X_2\oplus K Y_2\oplus \cdots \oplus K X_{m+1}\oplus K Y_{m+1}$.  \newline
  Now it is easy to see that 
 $\n$ is the $(m+1)$-th Heisenberg algebra over $K$.


\end{proof}

Recall that  $r_{1,\C}(\mathcal N)$  and   $W_{1, \C}$ are defined before Lemma \ref{crank1}.

\b{Le}\label{rankspan2}
    Let $\mathcal N$ be a Carnot algebra. Suppose $r_1(\mathcal N)\ge 2$  and
  $r_{1,\C}(\mathcal N)\le 1$.  Assume the action of  $\text{Aut}_g(\mathcal{N})$  on $V_1$  is irreducible.
 Then $V_1$   can be written as a direct sum of vector subspaces
  $V_1=U_1\oplus \cdots \oplus U_n$ such that
$[U_j, U_{j'}]=0$ for $j\not=j'$,  and
 for each $j$, $<U_j>$ is isomorphic
 to a complex Heisenberg algebra.
\end{Le}

\b{proof}
Lemma \ref{crank1} implies  $r_{1, \C}(\mathcal N)=1$. 
 We  claim   that  $W_{1, \C}=V_1\otimes \C$  holds.
 Suppose $W_{1, \C}\not=V_1\otimes \C$.
  Define  projections $\pi_1: W_{1, \C}\ra V_1$  and  $\pi_2: W_{1, \C}\ra V_1$  by $\pi_1(X_1+i X_2)=X_1$  and    $\pi_2(X_1+i X_2)=X_2$  respectively.     Clearly    $\pi_1$
  and $\pi_2$   commute    with  the action of
   $\text{Aut}_g(\mathcal{N})$. 
Since  $W_{1, \C}$   is invariant under the action of  $\text{Aut}_g(\mathcal N)$, 
  we see that  $\pi_1(W_{1, \C})$ is a subspace of $V_1$ invariant under the action of  $\text{Aut}_g(\mathcal{N})$.   Notice that $\pi_1(W_{1, \C})$ is non-trivial since $W_{1, \C}$ is non-trivial and is a complex subspace
of $V_1\otimes \C$  (so if $X_1+iX_2\in W_{1, \C}$, then $i(X_1+i X_2)\in W_{1, \C}$).   Since by our assumption
 the action of  $\text{Aut}_g(\mathcal{N})$  on $V_1$  is irreducible, we must have $\pi_1(W_{1, \C})=V_1$.
   Similarly the kernel $\ker (\pi_1)\subset i V_1$ is also invariant under
the action of  $\text{Aut}_g(\mathcal{N})$.  Hence $\ker (\pi_1)=\{0\}$ or $iV_1$.
     Since  we assumed   $W_{1, \C}\not=V_1\otimes \C$,   we must have
 $\ker (\pi_1)=0$.      A similar argument shows  that   $\pi_2: W_{1, \C}\ra V_1$ is onto  with trivial kernel. It follows that
  there is a  (real)  linear isomorphism $J: V_1\ra V_1$ such that $W_{1, \C}=\{X+i J(X):  X\in V_1\}$.
  Notice that, if $X=X_1+i X_2\in V_1\otimes \C$ has rank $1$, then its  complex     conjugate $\ol{X}=X_1-i X_2$ also   has  rank $1$.
 Since by definition $W_{1, \C}$ is spanned by rank one elements,    $\ol W_{1, \C}=\{X_1-iX_2:  X_1+iX_2\in W_{1, \C}\}$    is also  spanned by rank $1$ elements.  It follows that
$\ol W_{1, \C}\subset W_{1, \C}$.  Consequently,  $ W_{1, \C}=V_1\otimes\C $,     contradicting the facts that $\pi_1$ and $\pi_2$ have trivial kernels.
  The contradiction implies   $W_{1, \C}=V_1\otimes \C$.  

Now we apply Lemma   \ref{rankspan}
   to   $\n\otimes \C$ and 
conclude that $V_1\otimes \C$ is a  direct sum of   extended equivalence classes
  $V_1\otimes \C=W_1 \oplus \cdots \oplus W_m$,   and  $<W_j>$ is isomorphic to a complex Heisenberg algebra.
Note that the conjugate   $\ol W_j$ of $W_j$ is also an extended equivalence class, so it equals some $W_k$.    We claim that $\ol W_j\cap W_j=\{0\}$.
   Suppose $\ol W_j\cap W_j\not=\{0\}$.    Fix any
$0\not=X=X_1+i X_2 \in \ol W_j\cap W_j$,  where $X_1, X_2\in V_1$.
  Then   $\ol X=X_1-iX_2\in W_j$. Since $W_j$ is a complex linear subspace of $V_1\otimes \C$,   we have     $X_1, X_2\in W_j$.
 At least one of $X_1, X_2$ is nonzero.  We may assume $X_1\not=0$.
   Since $W_j$ is an extended  equivalence class,   every
nonzero element in $W_j$ has rank $1$.
    Hence there is some $Y\in V_1$ with $[X_1, Y]\not=0$,
   otherwise $[X_1, \n \otimes \C]=0$.
It follows that $[X_1, \n \otimes \C]=\C [X_1, Y]$.
Therefore
  $[X_1, \n]\subset \C [X_1, Y]\cap \n=\R [X_1, Y]$.
  This means $X_1\not=0$ has rank at most $1$, contradicting the assumption  that
  $r_1(\n)\ge 2$.   Hence $\ol W_j\cap W_j=\{0\}$.  Therefore
   $\ol W_j=W_k$ for some $k\not=j$.
 In particular $m=2n$ for some $n$.
We may relabel the $W_j$'s so that
  $W_{n+j}=\ol W_j$ for $1\le j\le n$.

Consider $W_j\oplus \ol W_{j}$ ($1\le j\le n$).
  Let $U_j=\pi_1(W_j)\subset V_1$.  Then      $W_j\oplus  \ol W_{j}=U_j\oplus i U_j$.  Hence $\dim_\C(W_j\oplus  \ol W_{j})=\dim_\R U_j$.
Notice that $V_1=U_1+\cdots + U_n$.
 Since $V_1\otimes \C=(W_1\oplus  \ol W_{1})\oplus \cdots\oplus  (W_n\oplus  \ol W_{n})$  and $\dim_\C(V_1\otimes \C)=\dim_\R V_1$, we obtain
  $\dim_\R V_1= \sum_j \dim_\R U_j$.  It follows that
  $V_1=U_1\oplus \cdots \oplus U_n$.   Lemma \ref{l21}  implies
 $[U_j, U_{j'}]=0$ for $j\not=j'$.
 It remains to show that $<U_j>$ is isomorphic to  a complex Heisenberg algebra.

 Recall that  $<W_j>$ is isomorphic to    ${\mathcal H}^k_\C$ for some $k \ge 1$.    Hence there   are  elements   $A_s, B_s\in   W_j$,  $1\le s\le k$, and $Z\not=0$ in the second layer of $<W_j>$ such that\newline
(1) $A_s, B_s$,  $1\le s\le k$,   
      form a basis of  $W_j$, and $Z$  spans the second layer of $<W_j>$; \newline
  (2)  $[A_s, B_s]=Z$  for all $s$,   and $[A_s, A_{s'}]=[B_s, B_{s'}]=[A_s, B_{s'}]=0$   for     $s\not=s'.$\newline
  By Lemma \ref{l21},
 $[W_j, W_{j'}]=0$ if  $j\not=j'$.  Since $\ol{W_j}=W_{j+n}$, we have
   $[A_s,  \ol A_s]=[A_s,  \ol B_s]=[\ol A_s, B_s]=[B_s, \ol B_s]=0$
 for all $s$,
  and  $[A, B]=0$ if   $s\not=s'$  and $A\in\{A_s, B_s, \ol A_s, \ol B_s\}$ and
 $B\in  \{A_{s'}, B_{s'}, \ol A_{s'}, \ol B_{s'}\}$.    Set
 $$X_s=\frac{1}{2}(A_s+\ol A_s), \;\;  \tilde X_s=\frac{1}{2i}(A_s-\ol A_s);$$
$$Y_s=\frac{1}{2}(B_s+\ol B_s), \;\;  \tilde Y_s=\frac{1}{2i}(B_s-\ol B_s);$$
$$ Z_1=\frac{1}{4}(Z+\ol  Z), \;\;  \tilde  Z_1=\frac{1}{4i}(Z-\ol Z).$$
      Since $W_j$ and   $\ol W_j$ are distinct extended equivalence classes, we have 
 $\C Z\cap \C \ol Z=\{0\}$. This implies  that $Z_1$ and $\tilde Z_1$ are linearly independent over $\R$.
Now it is easy to check the following:\newline
 (1) the following  bracket relations hold:
 $$[X_s, Y_s]=Z_1,\; [X_s, \tilde Y_s]= [\tilde X_s, Y_s]=\tilde Z_1,\; [\tilde X_s, \tilde Y_s]=-Z_1, [X_s, \tilde X_s]=[Y_s, \tilde Y_s]=0;$$
 $$[X, Y]=0 \;\;  \text{if}\; s\not=s',\; \text{and}\; X\in \{X_s, Y_s, \tilde X_s, \tilde Y_s\},    \;   Y\in \{X_{s'}, Y_{s'}, \tilde X_{s'}, \tilde Y_{s'}\}.$$
(2) $<U_j>$ is $2$-step,  $X_s, \tilde X_s, Y_s, \tilde Y_s$,  $1\le s\le k$,   form a basis of  
 $U_j$,  and
  $Z_1, \tilde Z_1$ form  a basis of  the second layer of  $<U_j>$.\newline
   It follows    that $<U_j>$ is isomorphic to the complex Heisenberg algebra
 $\mathcal H^k_\C$.

\end{proof}

Further information about the Carnot algebras in Lemma \ref{rankspan2}
 will be provided in Theorem \ref{nonrex}.  In fact, Theorem \ref{nonrex}
    gives  a characterization of  non-rigid Carnot   algebras   with
    irreducible  first layer.


Heisenberg product  algebras  and complex Heisenberg product  algebras
  are defined in the Introduction.    We use the notation from there.

\b{Le}\label{realhpg}
  Every Heisenberg product  algebra  $\n$  has an irreducible first layer and satisfies $r_1(\n)=1$.  In particular, it is non-rigid.
\end{Le}

\b{proof}
Let $P: \tilde \n\ra \n$ be the natural projection.   Then condition (1)
   in the definition of Heisenberg product algebra  implies that 
  the restriction $P|_{\tilde{\mathcal H}_j}$ is a graded isomorphism   from $\tilde{\mathcal H}_j$  onto  its image.   Set $L_j=P(\tilde L_j)$. Then  $\dim(L_j)=1$ and (2) implies $L_i\cap L_j=\{0\}$ for $i\not=j$.

 We first show that $r_1(\n)=1$.
Let $0\not=X\in V_1$. Then $X=P(\tilde X)$ for some $\tilde X\in \tilde V_1$.
   We claim that $\text{rank}(X)\ge 1$ always holds,  and $\text{rank}(X)=1$ if and only if
 $\tilde X\in  \tilde  V_{1,j}$ for some $j$.  The claim  implies  $r_1(\n)=1$ and so $\n$ is non-rigid by Theorem \ref{nonrch}.
 First assume   
$\tilde X\in \tilde V_{1,j}$  for some $j$.
Since $\text{rank}(\tilde X)=1$, we have $\text{rank}(X)\le 1$.    There is some $\tilde Y\in \tilde V_{1,j}$ such that $[\tilde X, \tilde Y]\not=0$.   Since
$P|_{\tilde{\mathcal H}_j}$ is a graded isomorphism, we have
$[X, P(\tilde Y)]=[P(\tilde X), P(\tilde Y)]=P([\tilde X, \tilde Y])\not=0$, which implies
 $\text{rank}(X)\ge 1$. Hence $\text{rank}(X)=1$.
  Next we assume   
$\tilde X\notin \cup_j\tilde V_{1,j}$.
  Write $\tilde X=\tilde X_{j_1}+\cdots +\tilde X_{j_s}$  with   $s\ge 2$  and 
$0\not=\tilde X_{j_i}\in \tilde V_{1,j_i}$,
  $1\le j_1< \cdots <  j_s\le n$.  
  There is some $\tilde Y_{j_i}\in  \tilde V_{1, j_i}$  satisfying
 $0\not=[\tilde X_{j_i}, \tilde Y_{j_i}]\in \tilde L_{j_i}$.    Then
  $0\not=[\tilde X, \tilde Y_{j_i}]=    [\tilde X_{j_i}, \tilde Y_{j_i}]  \in \tilde L_{j_i}$.
  We have $0\not=[X, P(\tilde Y_{j_i})]\in  L_j$.
 Since $L_i\cap L_j=\{0\}$ for $i\not=j$, we see that
$[X,  P(\tilde Y_{j_1})]\in L_{j_1}$ and
$[X,  P(\tilde Y_{j_2})]\in L_{j_2}$  are linearly independent. Hence $\text{rank}(X)\ge 2$.

 Next we show that $\n$  has irreducible first layer.
  Notice that,  each nonsingular linear transformation  of the second layer
  of the   Heisenberg algebra $\mathcal H^m_\R$,  which is just multiplication by a nonzero constant,     
extends to a  graded  automorphism of  $\mathcal H^m_\R$.  Hence 
  condition (4) in the definition of an   Heisenberg product algebra  implies each $B\in G_2$ extends to a
   graded  automorphism of $\tilde \n$.
  Now for each $B\in G_2$,  fix some  $\tilde B\in \text{Aut}_g(\tilde \n)$
 such that  $\tilde B|_{\tilde V_2}=B$.    Let $\tilde G=<\tilde B:  B\in G_2>\subset \text{Aut}_g(\tilde \n)$ be the group generated by all the $\tilde B$, $B\in  G_2$.     Notice that  the $\tilde G$-action on $\tilde V_2$ factors through   $G_2$.
 Therefore  condition  (3)  implies   that  $V$ is   $\tilde G$-invariant. 
 It follows that  each  $\tilde L\in  \tilde G$ induces a  graded isomorphism
 $L\in \text{Aut}_g(\n)$.    Let $G=\{L:  \tilde L\in \tilde G\}$. Then
 $G$ is a subgroup of $\text{Aut}_g(\n)$  and  the map  $\tilde G\ra G$,   $\tilde L\ra L$ is a homomorphism.     Condition (4) implies  that 
   $G$  permutes the set $\{P(\tilde V_{1,1}), \cdots,  P(\tilde V_{1,n})\}$ and acts transitively on the set.

  Let  $G_0\subset \text{Aut}_g(\h_\R^m)$ be the subgroup   consisting of graded automorphisms 
 that act trivially on the second layer.   Notice that $G_0$ acts   transitively on the set of nonzero vectors in
  the first layer (this follows  from the proof of  Lemma \ref{heisenberg}).  
  Since $G\subset  \text{Aut}_g(\n)$ acts on the set   $\{P(\tilde V_{1,1}), \cdots,  P(\tilde V_{1,n})\}$
 transitively,    we see that   $ \text{Aut}_g(\n)$      acts  transitively  on the nonzero vectors in 
  $\cup_j  P(\tilde V_{1,j})$.

  Now let  $U\subset V_1$ be a non-trivial subspace of   $V_1$ that is invariant under the action of $\text{Aut}_g(\n)$.    
       If  $U\cap   P(\tilde V_{1,j})\not=\{0\}$  for some $j$, then   $ P(\tilde V_{1,j})\subset U$ for all $j$ and so $U=V_1$.   Now assume $U\cap P(\tilde V_{1,j})=\{0\}$  for all $j$.
Fix an inner product on $V_1$  and   
    let $S\subset V_1$ 
 be the unit sphere   with respect to this inner product.  Set $K=S\cap (\cup_j P(\tilde V_{1,j}))$. Then
 $\delta:=d(U\cap S,   K)>0$.    Pick any $0\not=x\in U$.   Then  $x=P(\tilde x)$ for some $\tilde x\in \tilde V_1$.  
Write $\tilde x=\tilde x_{j_1}+\cdots + \tilde  x_{j_s}$ with  $s\ge 2$  and 
  $1\le j_1 <\cdots < j_s\le n$  and $0\not=\tilde x_{j_i}\in  \tilde V_{1, j_i}$.  
There exists   a sequence   
 $g_k\in G_0 $   such that    $g_k(\tilde x_{j_1})=k \tilde x_{j_1}$.  
    Define  $f_k\in   \text{Aut}_g(\tilde \n)$  by $f_k|_{\tilde{\h}_{j_1}}=g_k$ and 
$f_k|_{\tilde{\h}_j}=id$ for $j\not=j_1$.   Notice  that $f_k$ pointwise  fixes $\tilde V_2$.   Hence it  induces a graded  automorphism $F_k$ of $\n$.   
  We have $F_k(x)=k P(\tilde x_{j_1})+P(\tilde x_{j_2})+\cdots +P(\tilde x_{j_s})$.  
  It is clear that 
 $d(\R F_k(x)\cap S,   K)\ra 0$  as $k\ra \infty$. 
   Since $\R F_k(x)\subset U$, this contradicts 
 $\delta=d(U\cap S,   K)>0$.   Hence   the first layer of $\n$ is irreducible.

\end{proof}

  Notice  that, the map $\tau:  \h_\C^m=\C^{2m}  \oplus  \C\ra  \h_\C^m=\C^{2m}  \oplus  \C$,
  $$\tau(w_1, \cdots, w_{2m}, z)=(\bar w_1, \cdots, \bar w_{2m}, \bar z),$$ 
 is a graded automorphism of 
$ \h_\C^m=\C^{2m}  \oplus  \C$.

\b{Le}\label{complexhpg}
  Every complex  Heisenberg product  algebra  $\n$  has an irreducible first layer and satisfies $r_1(\n)=2$,
  $r_{1, \C}(\n)=1$.  In particular, it is non-rigid.
\end{Le}

\b{proof}
Let $P: \tilde \n\ra \n$ be the natural projection.   Then condition (1) implies the restriction $P|_{\tilde{\mathcal H}_j}$ is a graded isomorphism   from $\tilde{\mathcal H}_j$  onto  its image.   Set $L_j=P(\tilde L_j)$. Then  $\dim_\R(L_j)=2$.  
Recall that  we have $r_1(\h_C^m)=2$  and  $r_{1, \C}(\h_C^m)=1$  for   
    each complex Heisenberg algebra $\h_C^m$.  
Since $[P(\tilde  V_{1, i}),  P(\tilde  V_{1, j})]=0$ for $i\not=j$  and    
the restriction $P|_{\tilde{\mathcal H}_j}$ is a graded isomorphism   from $\tilde{\mathcal H}_j$  onto  its image, 
  it is easy to see that  
$r_1(\n)=2$,
  $r_{1, \C}(\n)=1$.  By Theorem \ref{nonrch},    $\n$ is non-rigid.

Next we show that $\n$  has irreducible first layer.
First we notice that,   if   we identify  $\C$ with the second layer of 
 $\h_\C^m=\C^{2m}  \oplus  \C$  and  
a  map $g: \C\ra \C$   has the form $g(z)=az$ or $g(z)=a \bar z$ for some $0\not=a\in \C$,   then $g$ extends to a graded isomorphism 
  $\h_\C^m=\C^{2m}  \oplus  \C\ra  \h_\C^m=\C^{2m}  \oplus  \C$.  
It follows that every $g_2\in  G_2$ extends to a  graded  automorphism of $\tilde \n$.  
   Then  the arguments in the last 3 paragraphs of the proof   of  Lemma \ref{realhpg}  show that $\n$ has an irreducible first layer.


\end{proof}

Finally we are ready to characterize  the non-rigid Carnot   algebras with irreducible first layer. 

\b{Th}\label{nonrex}
Let $\mathcal{N}$  be a non-rigid   Carnot algebra.
 Then
 $V_1$  is  irreducible  if and only if  exactly one of the following
 happens:    \newline
(1) $\n$ is abelian;\newline
  (2) $r_1(\mathcal{N})=1$   and   $\mathcal{N}$ is     a    Heisenberg product  algebra;\newline
(3) $r_1(\mathcal{N})=2$,    $r_{1, \C}(\mathcal{N})=1$  and   $\mathcal{N}$ is    a     complex Heisenberg  product 
  algebra.

\end{Th}

\b{proof}
  An abelian Carnot algebra   clearly is non-rigid and  has irreducible first layer.  Lemmas  \ref{realhpg}
  and \ref{complexhpg}  imply Heisenberg product  algebras  and  complex Heisenberg  product 
  algebras   are non-rigid and  have   irreducible first layer.

  Conversely,  let
     $\mathcal{N}$  be a non-abelian non-rigid   Carnot algebra  with  irreducible first layer.
  We need to show that $\n$   satisfies either (2) or (3). 
  Lemma \ref{rank0}   and 
 Theorem \ref{nonrch}  imply     that  we have one of the following
   two cases: (1) $r_1(\mathcal{N})=1$;  (2)
 $r_1(\mathcal{N})\ge 2$   and   $r_{1, \C}( \mathcal{N})\le 1$.

\noindent
{\bf{Case (1)  $r_1(\mathcal{N})=1$.}}

By  Lemma   \ref{rankspan}  and Lemma \ref{l21},
 $V_1$ can be  written as a direct sum of vector subspaces
  $V_1=U_1\oplus \cdots \oplus U_n$, where each $U_j$ is an
extended   equivalence  class,  $[U_i, U_j]=0$ for $i\not=j$; furthermore,
   $U_j$ is maximal among linear subspaces consisting of $0$ and
     rank $1$ elements,  and $<U_j>$ is isomorphic to an Heisenberg algebra.

Since  $U_j$ is maximal among linear subspaces consisting of $0$ and
     rank $1$ elements,
     the action of $\text{Aut}_g(\mathcal{N})$ on $V_1$  permutes the subspaces $U_1, \cdots,  U_n$.  Let $U$ be the direct sum of all those $U_j$ in the orbit of $U_1$.
  Then $U$ is invariant under the action of
      $\text{Aut}_g(\mathcal{N})$.
  Since  $V_1$ is irreducible,   we must have
 $U=V_1$.   Hence
$\text{Aut}_g(\mathcal{N})$    acts on the set  $\{U_1, \cdots,  U_n\}$  transitively.
  It follows that
    all the subalgebras $<U_j>$ are isomorphic.
 Suppose  they are isomorphic to  ${\mathcal H}_\R^m$  
for some $m\ge 1$.

Let $\tilde \n=\tilde {\mathcal H}_1\oplus \cdots \oplus\tilde  {\mathcal H}_n$,
  where   each $\tilde {\mathcal H}_j={\mathcal H}_\R^m$.   
Let $\tilde V_{1,j}$ and  $\tilde L_j$ respectively  be the first and second layers of   $\tilde {\mathcal H}_j$.  Let $\tilde V_1$ and $\tilde V_2$ respectively  be the first and second layers of  $\tilde \n$.
  Then     $\tilde V_1=\tilde V_{1,1}\oplus \cdots \oplus \tilde V_{1,n}$
  and $\tilde V_2=\tilde L_1\oplus \cdots \oplus \tilde L_n$.
  For each $j$, let $f_j: \tilde {\mathcal H}_j \ra <U_j>$ be a fixed
 graded isomorphism. We also use $f_j$ to denote the composition  of
$f_j:  \tilde {\mathcal H}_j \ra <U_j>$ with the inclusion $<U_j>\subset \n$.
 Now define $P: \tilde \n\ra \n$ by $P(x_1, \cdots, x_n)=\sum_j f_j(x_j)$.
Using the facts that $[U_i, U_j]=0$ for $i\not=j$  and that
$V_1=U_1\oplus \cdots \oplus U_n$,   it is easy to  check that $P$
is a surjective  graded homomorphism, and    $V:=\ker P\subset \tilde V_2$.
   It follows that
 $\mathcal{N}\cong  \tilde \n/V$.

   As  $P|_{\tilde {\mathcal H}_j} =f_j$ is an isomorphism,     we have  $V\cap \tilde L_j=\{0\}$.
Set   $L_j=f_j(\tilde  L_j)$.    Note that $L_j$  is the second layer of
 $<U_j>$.
  Since $U_i$ and $U_j$ are distinct extended equivalence classes for $i\not=j$,    we have $L_i\cap L_j=\{0\}$.
  This implies  $(V+\tilde L_i)\cap (V+\tilde L_j)=V$ for $i\not=j$.


   Since  $ \text{Aut}_g(\mathcal{N})$  permutes the   subalgebras $<U_j>$
     and
 $P|_{{\tilde {\mathcal H}}_j}=f_j$ is an isomorphism  from   ${\tilde {\mathcal H}}_j$
onto $<U_j>$,     for each
    $A\in \text{Aut}_g(\mathcal{N})$,
 there is a unique lift
    $\tilde A_j: {\tilde {\mathcal H}}_j\ra \tilde \n$
for $A|_{<U_j>}: <U_j>\ra \n $.    In fact, if $A(<U_j>)=<U_{\sigma(j)}>$ for some $\sigma(j)$,  then
 $\tilde A_j=f^{-1}_{\sigma(j)}\circ A\circ f_j$.
  Define $\tilde A:  \tilde \n\ra \tilde \n$ by $\tilde A(x_1,  \cdots, x_n)=\sum_j \tilde A_j(x_j)$.
  Then it is easy to see that $\tilde A$ is a graded isomorphism and is a lift of $A$(i.e.,
 $P\circ \tilde A= A\circ P$).    Furthermore,     it is also easy to check that $\tilde A$ is the only graded isomorphism that lifts $A$.
It follows that  the map   $h:  \text{Aut}_g(\mathcal{N})\ra
\text{Aut}_g(\tilde \n)$,
  $h(A)=\tilde A$, is an injective homomorphism.
  Set $G=h(\text{Aut}_g(\mathcal{N}))$.
  Then  $V$ is $G$-invariant.
   Since  $\text{Aut}_g(\mathcal{N})$ acts transitively on the set $\{U_1, \cdots, U_n\}$,
  $G$ acts transitively on the set  $\{\tilde {\mathcal H}_1, \cdots, \tilde {\mathcal H}_n\}$.

    Let
$R:\text{Aut}_g(\tilde {\mathcal{N}})\ra GL(\tilde V_2)$ be defined by
$R(\tilde A)=\tilde A|_{\tilde V_2}$.  Set $G_2=R(G)$. Clearly the $G$ action  on
$\tilde V_2$ factors through $G_2$.  Hence $V\subset \tilde V_2$ is
$G_2$-invariant and    $G_2$ acts transitively on the set $\{\tilde L_1,
\cdots,  \tilde L_n\}$ of \lq\lq coordinate axes" of $\tilde V_2$.


\noindent {\bf{Case (2)  $r_1(\mathcal{N})\ge 2$   and   $r_{1, \C}(
\mathcal{N})\le 1$}}.

The proof is similar to that of case (1). We will mainly indicate
the difference.

Lemma \ref{rankspan2}  implies
   $V_1$   can be written as a direct
sum of vector subspaces
  $V_1=U_1\oplus \cdots \oplus U_n$ such that
$[U_j, U_{j'}]=0$ for $j\not=j'$,  and
 for each $j$, $<U_j>$ is isomorphic
 to a complex Heisenberg algebra.

Recall that  (see the proof of Lemma \ref{rankspan2})
 (a)
$U_j=\pi_1(W_j)$, where $W_j$ is an extended equivalence class  of
$V_1\otimes \C$; (b)  each extended equivalence class  is a maximal
linear subspace (over $\C$) of $V_1\otimes \C$
  consisting of ($0$ and)  rank $1$ elements.
  Hence  $\text{Aut}_g(\mathcal{N})$  permutes
    the $W_j$'s. 
Notice that for any $A\in \text{Aut}_g(\mathcal{N})$, we have
 $A\circ \pi_1=\pi_1\circ  A$.    
  Suppose $A(W_j)=W_{\sigma(j)}$  for   some $\sigma(j)$.   
   Then    $A(U_j)=U_{\sigma(j)}$.   So
 $\text{Aut}_g(\mathcal{N})$  permutes
    the $U_j$'s.
      Since  the action of 
$\text{Aut}_g(\mathcal{N})$   on $V_1$ is irreducible,  we see that 
$\text{Aut}_g(\mathcal{N})$ 
 acts transitively on the set
$\{U_1, \cdots, U_n\}$.  It follows that
    all the subalgebras $<U_j>$ are isomorphic.
 Suppose  they are isomorphic to  ${\mathcal H}^m_\C$  for some $m\ge 1$.

Let $\tilde \n=\tilde {\mathcal H}_1\oplus \cdots \oplus \tilde {\mathcal H}_n$,
  where   each $\tilde {\mathcal H}_j={\mathcal H}^m_\C$.   
Let $\tilde V_{1,j}$ and  $\tilde L_j$ respectively  be the first and second layers of   $\tilde {\mathcal H}_j$.  Let $\tilde V_1$ and $\tilde V_2$ respectively  be the first and second layers of  $\tilde \n$.
  Then     $\tilde V_1=\tilde V_{1,1}\oplus \cdots \oplus \tilde V_{1,n}$
  and $\tilde V_2=\tilde L_1\oplus \cdots \oplus \tilde L_n$.
  For each $j$, let $f_j: \tilde {\mathcal H}_j \ra <U_j>$ be a fixed
 graded isomorphism. 
As in Case (1) we define a  surjective  graded homomorphism $P:
\tilde \n\ra \n$.  Let $V=\ker P$.   Then  $V\subset \tilde V_2$  and
 $\mathcal{N}\cong  \tilde \n/V$.
Since $P|_{\tilde {\mathcal H}_j}=f_j$ is an
isomorphism,  we see that  $\tilde L_j\cap V=\{0\}$.

As in the proof of Case (1), we construct  groups  $G\subset \text{Aut}_g(\tilde \n)$ and 
 $G_2\subset GL(\tilde V_2)$. The arguments there  show that $V$ is $G_2$-invariant,  
 $G_2$ permutes the set $\{\tilde L_1, \cdots,  \tilde L_n\}$ and acts transitively on the set. 
Recall that each $g_2\in G_2$ is the restriction  $g|_{\tilde V_2}$   of some $g\in G\subset \text{Aut}_g(\tilde \n)$. 
 Notice that $g$ permutes the $\tilde\h_j$'s and for each $j$,     $g|_{\tilde\h_j}: \tilde\h_j\ra \tilde\h_{\sigma(j)}$ is an isomorphism of the complex Heisenberg algebra.   By  \cite{S},   for every graded automorphism $A$ 
of the   complex Heisenberg algebra $\h_\C ^m=\C^{2m}\oplus \C$,  the map $B:=A|_{\{0\}\oplus \C}:  \C\ra \C$ has the form 
 $B(z)=az$ or $B(z)=a \bar z$ for some $0\not= a\in \C$.  Hence $V$ and $G_2$ also satisfy condition (4) in the definition of the complex Heisenberg product algebra.

\end{proof}

\section{Rigidity of quasisymmetric maps}\label{highstep}

  The goal of this Section  is to show that  non-rigid  Carnot   algebras  with  reducible first layer are 
 quasisymmetrically rigid.  
We  consider two cases depending on whether the invariant subspace is abelian.

\subsection{When the invariant subspace  is non-abelian}\label{nonabelian}

In this Subsection we consider the case when 
  the invariant subspace  is non-abelian.

\b{Le}\label{r=1}
  Let $\n$ be a $r$-step Carnot algebra with $r\ge 3$.  If 
   $r_1(\mathcal N)=1$   and   
${\mathcal W}:=<W_1>$ is not abelian, then  $\n$ is
  quasisymmetrically rigid.

\end{Le}

\b{proof}  
By   Corollary \ref{highrank},
    $W_1$ is a non-trivial proper subspace of $V_1$  invariant under the action of   $\text{Aut}_g(\mathcal N)$.
  Since  $r_1(\mathcal N)=1$  and $W_1$ is spanned by rank $1$ elements,
 we have $r_1(\mathcal W)\le 1$.  First assume
$r_1(\mathcal W)=0$.  Lemma \ref{rank0.1}  implies $\mathcal W=\R^m\oplus \mathcal W'$ is a direct sum of an Euclidean algebra and  another Carnot algebra. 
  Since   $\mathcal  W$ is not abelian,  $\mathcal W'$ is not Euclidean. 
   Then the main result in \cite{X3}   implies  $\mathcal W$ is
   quasisymmetrically rigid.
  Now Theorem \ref{2step3}  implies $\n$ is quasisymmetrically rigid.

  Next we assume
$r_1(\mathcal W)=1$.
 We claim that $[X, W_1]\subset \mathcal W$  for any $X\in V_1$.  The Lemma then follows from   Theorem \ref{2step2}.
  Indeed,
  by definition $W_1$ has a vector space basis $\{v_1, \cdots, v_m\}$ such that each $v_i$ has rank $1$ in the algebra $\mathcal N$. Since
$r_1(\mathcal W)=1$, $v_i$ also   has  rank $1$ in the  algebra
$\mathcal W$.  Hence there is some $w_i\in W_1$ such that
 $[v_i, w_i]\not=0$.
  Since $v_i$ has rank $1$ in $\mathcal N$, for any $X\in V_1$, we have
 $[v_i, X]\in \R [v_i, w_i]\subset [W_1, W_1]$.  Since this holds for all $i$, we see that $[X, W_1]\subset [W_1, W_1]$ for all $X\in V_1$.  Hence  the claim holds.

\end{proof}

The subspace $\hat W_1$  of $V_1$    is defined in Corollary \ref{crank1cor}.

\b{Le}\label{r=2}
Let $\n$ be a $r$-step Carnot algebra with $r\ge 3$.
 If
 $r_1(\mathcal N)\ge 2$,
$r_{1,\C}(\mathcal N)\le 1$,
 and $\hat {\mathcal W}:=<\hat W_1>$ is not abelian, then  $\n$ is
  quasisymmetrically rigid.

\end{Le}

\b{proof}  
Lemma \ref{crank1}
  and Lemma \ref{l4.12} imply that  $r_1(\mathcal N)= 2$  and
$r_{1,\C}(\mathcal N)= 1$.
Corollary \ref{crank1cor}  implies  $\hat W_1$
 is a non-trivial proper subspace of $V_1$  invariant under the action of   $\text{Aut}_g(\mathcal N)$.   
In the proof of Lemma \ref{l4.12}, we   showed that if $X=X_1+iX_2\in V_1\otimes \C$
 satisfies $\text{rank}(X)=1$, then $\text{rank}(X_1)\le 2$.  Since $r_1(\n)=2$, we have either $X_1=0$ or $\text{rank}(X_1)=2$.   Since $W_{1, \C}$ is spanned by rank one elements in $V_1\otimes \C$, we see that $\hat W_1$   is spanned by  rank two  elements
 in  $\n$.  It follows that $r_1(\hat{\mathcal W})\le 2$.

 {\bf{Claim}} $r_1(\hat{\mathcal W})=0$   or $r_1(\hat{\mathcal W})=2$.

We shall first
 assume the claim and finish the proof of the Lemma, then  prove the claim at the end.
    First assume
$r_1(\hat{\mathcal W})=0$.
  Since   $\hat{\mathcal  W}$ is not abelian,  the main result in \cite{X3}  implies  $\hat{\mathcal W}$ is
   quasisymmetrically rigid.
  Now Theorem \ref{2step3}  implies $\n$ is quasisymmetrically rigid.

  Next we assume
$r_1(\hat{\mathcal W})=2$.
 We claim that $[X, \hat W_1]\subset \hat{\mathcal W}$  for any $X\in V_1$.  The Lemma then follows from   Theorem \ref{2step2}.
  Indeed,  as observed above,
 $\hat W_1$ has a vector space basis $\{v_1, \cdots, v_m\}$ such that each $v_i$ has rank $2$ in the algebra $\mathcal N$. Since
$r_1(\hat{\mathcal W})=2$, $v_i$ also   has  rank $2$ in the  algebra
$\hat{\mathcal W}$.   It follows that $[v_i,  \n]=[v_i,  \hat{\mathcal W}]\subset \hat{\mathcal W}$.
 Since this holds for all $i$, we see that $[\hat W_1, \n]\subset \hat{\mathcal W} $.
   Hence  the claim holds.

 Now  we  prove the claim. 
We assume
$r_1(\hat{\mathcal W})\ge 1$  and need to show
$r_1(\hat{\mathcal W})=2$.
  Note that  if $A$ is an
extended equivalence class   in $V_1\otimes \C$, then so is its complex  conjugate  $\bar A$.
Furthermore,  $A\oplus \bar  A=U\oplus iU$, where
 $$U=\pi_1(A)=\{X\in V_1:  X+iY\in A \;\text{for some}\;  Y\in V_1\}.$$
Let $0\not=X_1\in U$.  Then  $X_1+iX_2\in A$ for some $X_2\in V_1$.
In the proof of Lemma \ref{l4.12} we  showed
 that   $\text{rank}(X_1)\le 2$.  Our assumption $r_1(\n)\ge 2$ now implies
$\text{rank}(X_1)=2$.
     Since $\text{rank}(X_1+iX_2)=1$, we have
 $[X_1+iX_2,  V_1]\not=0$, otherwise
$[X_1+iX_2,  \n\otimes \C]=0$.   Let $Y_0\in V_1$ be an arbitrary element such that
 $[X_1+iX_2, Y_0]\not=0$.
Then  $[X_1+iX_2,   \n]\subset [X_1+iX_2, \n\otimes \C]=\C [X_1+iX_2, Y_0]$,
 which implies $[X_1, \n ]\subset \R [X_1, Y_0]+\R [X_2, Y_0]$.
Since $\text{rank}(X_1)=2$, we see that $[X_1, Y_0]$ and $[X_2, Y_0]$ are linearly independent
  (over $\R$).

Notice that $W_{1, \C}$ is spanned by extended equivalence classes in $V_1\otimes \C$.
There are extended equivalence classes  $A_1, \cdots, A_n$  such that  
$W_{1, \C}=(A_1+\bar A_1)+\cdots +(A_n+\bar A_n)$  and 
$A_i\not=A_j, \bar A_j$ 
 for $i\not=j$.    Let  $U_j=\pi_1(A_j)$. Then  
$A_j\oplus \bar{A_j}=U_j\oplus iU_j$  and  
 $\hat{ W}_1=U_1+\cdots + U_n$.
 Lemma \ref{l21}
implies  $[U_i, U_j]=0$   for  $i\not=j$.
   Since  $r_1(\hat{\mathcal W})\ge 1$,  
    for any $0\not=Y_0\in U_j$,  there must exist some $X_1\in U_j$ such that $[X_1, Y_0]\not=0$.
The preceding  paragraph now implies  $[Y_0, U_j]=[Y_0, \n]$ has dimension $2$.
  Now any $0\not=Y\in \hat{W}_1$ can be written as  $Y=Y_{j_1}+\cdots  + Y_{j_s}$,
  where $1\le j_1< \cdots< j_s\le n$,  and $0\not=Y_{j_i}\in U_{j_i}$.
 It follows that
  $[Y, \hat{W}_1]\supset [Y, U_{j_1}]=[Y_{j_1},  U_{j_1}]$ has dimension  at least $2$.
  Hence $Y$ has rank at least $2$ in $\hat{\mathcal W}$.

\end{proof}

\subsection{When invariant subspace  is    abelian}\label{abelian}

In this subsection we   first treat the case when the invariant subspace  is    abelian, then finish the proof  of Theorem \ref{main}.

The following lemma follows from the   BCH formula (see the 
proof of Lemma 2.1 in \cite{X1}).

\b{Le}\label{cbhf}
There are universal constants $c_j\in \Q$ for   $j\ge 1$   with the following property:
  $c_1=1/2$,  and  
  if $\n$ is $r$-step  and $Y\in \n$ satisfies $[Y,  V_i]=0$ for all $i\ge 2$, then
     for any  $X\in \n$ we have:
\b{equation}\label{e2.3.10}
X*Y=X+Y+\sum_{j=1}^{r-1} c_j(\text{ad} \,X)^j \, Y
\end{equation}
  and
\b{equation}\label{e2.3.11}
Y*X=Y+X+\sum_{j=1}^{r-1} (-1)^j c_j(\text{ad} \,X)^j \, Y.
\end{equation}

\end{Le}

\b{Le}\label{l7.11} 
   Let $W\subset V_1$ be a linear subspace satisfying $[W, W]=0$ and $[W, V_i]=0$ for all $i\ge 2$. 
Let $\tilde W\subset V_1$ be a subspace complementary to $W$.  
  Then for any $u\in W$  and any 
$
x=x_1+\sum_{j=1}^r \tilde x_j$ ($x_1\in W$, $\tilde x_1\in \tilde
 W$,  $\tilde x_j\in V_j$ for $j\ge 2$),   we have
   $(ad\;  x)^i u=(ad\; \tilde x_1)^i u$ for all $i\ge 1$.
\end{Le}

\b{proof} We first show that for any integer $k\ge 0$,  any $j\ge 2$ and any $y_j\in V_j$ we have
  \b{equation}\label{e7.70}
  [y_j, (ad\; \tilde x_1)^ku]=0.
  \end{equation}
   We induct on $k$. The base case $k=0$ follows from the assumption 
that $[W, V_j]=0$  for $j\ge 2$.   Assume the claim holds
for $k$.  Now the Jacobi
    identity implies
 $$[y_j, (ad\; \tilde x_1)^{k+1} u]=[y_j,  [\tilde x_1,(ad\; \tilde
 x_1)^k u ]]=[[y_j, \tilde x_1], (ad\; \tilde x_1)^k  u]+[\tilde x_1,
 [y_j, (ad\; \tilde x_1)^k  u]].
 $$
    Notice that $[y_j, \tilde x_1]\in V_{j+1}$. The induction
    hypothesis now implies   
    $[[y_j, \tilde x_1], (ad\; \tilde x_1)^k  u]=0$ and $
 [y_j, (ad\; \tilde x_1)^k  u]=0$.

Next we prove $(ad\; \tilde x)^i u=(ad\; \tilde x_1)^i u$
 by inducting on $i$.
  The base case $i=1$ follows from the assumptions $[W, W]=0$ and $[W, V_j]=0$ for $j\ge 2$:
   $$(ad\; \tilde x) u=[\tilde x, u]=[x_1,
  u]+[\tilde x_1, u]+\sum_{j=2}^{r}[\tilde x_j, u]=0+ [\tilde x_1, u] +0=(ad\; \tilde x_1)u. $$
    Now assume  $(ad\; \tilde x)^i u=(ad\; \tilde x_1)^i u$  holds
    for $i$.    Then   the induction hypothesis, the condition $[W,
    V_j]=0$ for $j\ge 2$ and (\ref{e7.70}) imply
\b{align*}
(ad\; \tilde x)^{i+1} u=[\tilde x, (ad\; \tilde x)^i u]&=[\tilde x, (ad\; \tilde x_1)^i u]\\
 & =[x_1, (ad\; \tilde x_1)^i u]+[\tilde x_1, (ad\; \tilde x_1)^i u]+\sum_{j=2}^{r}[\tilde x_j,(ad\; \tilde x_1)^i u]\\
  &=0+(ad\; \tilde x_1)^{i+1} u+0=(ad\; \tilde x_1)^{i+1} u.
  \end{align*}

\end{proof}

\b{Th}\label{abelianhigh}  
Let  $\n=V_1\oplus \cdots\oplus V_r$  be  a  Carnot algebra,   $W\subset V_1$  a
non-trivial proper subspace
   invariant under  the action of $\text{Aut}_g(\mathcal{N})$.  If     $[W, \, W]=0$  and
  $[W, \, V_i]=0$ for all $i\ge 2$,  then   $\n$ is quasisymmetrically rigid. 
\end{Th}

Let $\tilde W\subset V_1$ be a subspace complementary to $W$.
We fix   an  inner product  $|\cdot|$   on  each of $W, \tilde W, V_j$, $j\ge 2$.  
Then there exists some constant $A\ge 1$ such that
 \b{equation}\label{e7.10}
|[v, w]|\le A \cdot |v|\cdot |w|
\end{equation}
   for all $v, w\in W\cup\tilde W\bigcup\cup_{j=2}^r V_j$.

Define a norm on
 $\n$   by
$$||x||=|x_1|+\sum_{j=1}^r |\tilde x_j|^{\frac{1}{j}}$$
   for  $ x=x_1+\tilde x_1+\cdots +\tilde x_r\in \n$  with $x_1\in W,  \tilde x_1\in \tilde W,  \tilde x_j\in V_j$     for 
   $ j\ge 2.$

We recall that,  if $F: X\ra Y$ is $\eta$-quasisymmetric, then $F^{-1}: Y\ra X$ is
  $\eta_1$-quasisymmetric  with $\eta_1(t)=(\eta^{-1}({t}^{-1}))^{-1}$.
  Without loss of generality we may assume $\eta(1)\ge 1$. It follows that  $\eta_1(1)\ge 1$.
 These inequalities will be used  implicitly.

The condition $[W, W]=0$ implies $W$ is a Lie subalgebra of $\n$.  We shall abuse notation and
  also   denote by $W$ the   connected Lie subgroup  of $N$   with Lie algebra $W\subset \n$. Here $N$ denotes the Carnot group with Lie algebra $\n$.

Let $d$ be the homogeneous distance on $N$ associated with the above defined norm on $\n$. That is,
  $d(x,y)=||(-x)*y||$  for $x,y\in \n=N$.   Let  $F:  (N, d)  \ra (N, d)$
  be an $\eta$-quasisymmetric map for some $\eta$. 
  Since $W\subset \n$ is an invariant subspace,  Proposition \ref{p3} 
    implies that $F$ sends   left cosets of $W$ to left   cosets of $W$.

Recall that $l_F$ and $L_F$ are defined in Section \ref{maps}.

\b{Le}\label{key}
Let $L$  be   a left coset of $W$.
   Suppose $p,q\in L$ are such that $l_F(p)>  C_1\cdot  L_F(q)$  with
$C_1= 200(r+1)A^{r-1}C \eta_1(1)$,  where $C=1+\max_{1\le j\le r-1} |c_j|$  and
  the  $c_j$'s    are the constants in Lemma \ref{cbhf}.
   Write $q=p*u$ with $u\in W$.   Then for any $s=p*tu$  with $|t|\ge 1$,    we   have
 $$ L_F(s)\le \frac{2(r+1)CA^{r-1} (\eta_1(1))^2}{|t|^{\frac{1}{r}}}\cdot
 L_F(q).$$


\end{Le}

\b{proof}  
  The proof is a modification of  the  proof of Lemma 5.1, \cite{X2}.

Fix some $e\in V_2$  with   $|e|=1$.
Set   $p'=F(p)$,   $q'=F(q)$   and  $L'=F(L)=p'*W$.
Let $\{r_j\}$ be an arbitrary sequence of positive reals such that $r_j\ra 0$.
Denote $\bar p_j=p'*r^2_j e$ and $\bar q_j=q'*r^2_j e$.  Note $p', q'\in  L'$. So
 $q'=p'*u'$ for some $u'\in W$.  Since $[W, V_2]=0$, we have $u'*r_j^2e=r_j^2e*u'$. 
It follows that 
 $$\bar q_j=p'*u'*r^2_j e=p'*r^2_j e*u'=\bar p_j*u'\in \bar p_j*W. $$ 
  Hence  $\bar p_j$ and $\bar q_j$ lie on the same left coset of $W$. 
 Let $p_j$  and   $ q_j$   
be points on $ F^{-1}(\bar p_j*W)$ nearest to $p$ and $q$ respectively.

Since  $l_F(p)>  C_1\cdot  L_F(q)$,  
  the  first   part of  the  proof of Lemma 5.1, \cite{X2}   yields 
\b{equation}\label{e3.1}
   {d(p, p_j)}\le \frac{1}{101(r+1)CA^{r-1}}d(q, q_j)
\end{equation}
  for all sufficiently large $j$.

Next we shall look at $d(p, p_j)$  and $d(q, q_j)$.

Since $q, p$ lie on the same left  coset,  we  can write
  $q=p*u$   for some $u\in W$.
  Similarly,   $q_j=p_j*w$  for some $w\in W$.
  Let $o_j=p_j*w'$ ($w'\in W$)  be an arbitrary point on the left coset $p_j*W$.
  Also
   write $p_j=p*(x_1+ \tilde x)$,   where
 $x_1\in W$  and
 $\tilde x=\tilde{x}_1+ \tilde{x}_2+\cdots +\tilde x_r$
  with  $\tilde x_1\in \tilde W$ and  $\tilde{x}_i\in  V_i$  for $i\ge 2$.   Although the  $x_1$   and
     $\tilde{x}_i$'s      depend on
 $r_j$,  we shall surpress the dependence to simplify the notation.


Now   we calculate $d(p, p_j)$  and $d(q, o_j)$.
  Notice
$$(-p)*p_j= (-p)*p*(x_1+ \tilde{x})=x_1+ \tilde{x}_1+ \tilde{x}_2+\cdots +\tilde x_r.$$
So
  \b{equation}\label{e3.2}
d(p, p_j)=|| (-p)*p_j||=|x_1|+\sum_{j=1}^r|\tilde x_j|^{\frac{1}{j}}.
  \end{equation}
Using Lemma \ref{cbhf},  the conditions $[W, W]=0$, $[W, V_i]=0$ for $i\ge 2$,  and  Lemma \ref{l7.11},    we obtain:
\b{align*}
  &(-q)*o_j  \\
& =(-q)*p_j*w'\\
&=(-u)*(-p)*p*(x_1+ \tilde{x})*w'\\
&=(-u)*(x_1+ \tilde{x})*w'\\
&=\left((x_1-u)+\tilde x-\sum_{j=1}^{r-1}(-1)^j c_j(ad \,(x_1+\tilde x))^j u\right)*w'\\
&=\left((x_1-u)+\tilde x-\sum_{j=1}^{r-1}(-1)^j c_j(ad \,\tilde x_1)^j u\right)*w'\\
&=(x_1-u+w')+\tilde x-\sum_{j=1}^{r-1}(-1)^j c_j(ad \,\tilde x_1)^j u
+\sum_{j=1}^{r-1}c_j(ad \,\tilde x_1)^j
w'\\
&=(x_1-u+w')+\tilde x_1+\sum_{j=1}^{r-1}\{\tilde x_{j+1}+c_j(ad \,\tilde
x_1)^j(w'-(-1)^ju)\}.
\end{align*}
  Hence
 \b{equation}\label{e3.5}
d(q, o_j)=|x_1-u+w'|+|\tilde x_1|+\sum_{j=1}^{r-1}\left|\tilde x_{j+1}+c_j(ad
\,\tilde x_1)^j(w'-(-1)^ju)\right|^{\frac{1}{j+1}}.
\end{equation}

Since $q_j$ is a point on $p_j*W$ nearest to $q$, we have
   $d(q_j, q)\le d(o_j, q)$ for any $o_j\in p_j*W$.   By (\ref{e3.1}), we have
 $d(p, p_j)\le d(q, o_j)/(101(r+1)C A^{r-1})$.
  In particular, this inequality holds for the point $o_j=p_j*u$.
Set $a=d(q, p_j*u)$.
Notice that
 \b{equation}\label{e7.19}
a=d(q, p_j*u)=|x_1|+|\tilde x_1| +\sum_{j: \,\text{even}}|\tilde
x_{j+1}|^{\frac{1}{j+1}}+\sum_{j:\, \text{odd}}|\tilde
x_{j+1}+2c_j(ad \,\tilde x_1)^j u|^{\frac{1}{j+1}},
\end{equation}
 where $j$ in the above sums varies from $1$ to $r-1$.
 Since $d(p, p_j)\le   a/(101(r+1)C A^{r-1})$,
  from (\ref{e3.2}) we get
\b{equation}\label{e7.20}
|x_1|\le a/(101(r+1)C A^{r-1})  \;\;\text{and}\;\;
  |\tilde x_j|^{\frac{1}{j}}\le a/(101(r+1)C A^{r-1})
\end{equation}
   for all $1\le j\le r$.
  Now   (\ref{e7.19})  and  (\ref{e7.20})  imply  
  \b{equation}\label{e7.21}
|\tilde x_{j_0+1}+2c_{j_0}(ad \,\tilde x_1)^{j_0}
u|^{\frac{1}{j_0+1}}\ge a/(r+1)
  \end{equation}
 for some  odd $j_0$,   $1\le j_0\le r-1$.  The triangle inequality
 and (\ref{e7.20}), (\ref{e7.21})  imply
  \b{equation}\label{e7.22}
  |2c_{j_0}(ad \,\tilde x_1)^{j_0}u|\ge
  \left(1-\frac{1}{(101C A^{r-1})^{j_0+1}}\right)\left(\frac{a}{r+1}\right)^{j_0+1}.
  \end{equation}

 Next we consider a point $s\in p*W$ of the form $s=p*tu$, where   $|t|\ge 1$.
  We claim that
\b{equation}\label{e3.6.5} d(s,\,  p_j*W)\ge
\frac{{|t|^{\frac{1}{j_0+1}}}}{2(r+1)CA^{r-1}} \cdot d(q, q_j).
\end{equation}
In the formula    (\ref{e3.5})  we replace $u$ with $tu$  and $q$
with $s$   to  obtain:
 \b{equation}\label{e3.9}
d(s, o_j)=|x_1-tu+w'|+|\tilde x_1|+\sum_{j=1}^{r-1}\left|\tilde
x_{j+1}+c_j(ad \,\tilde x_1)^j(w'-(-1)^j tu)\right|^{\frac{1}{j+1}}.
\end{equation}

Set $v=w'-tu$.
We divide points on $p_j*W$ into two types,  depending on whether
 $|v|\ge  |t|^{\frac{1}{j_0+1}}\cdot \frac{a}{(r+1)C A^{r-1}}$.   If   $|v|\ge  |t|^{\frac{1}{j_0+1}}\cdot \frac{a}{(r+1)CA^{r-1}}$, then
   (\ref{e3.9}) and
  (\ref{e7.20}) imply
\b{align}\label{e7.35}
 d(s, o_j)\ge |x_1-tu+w'|=|x_1+v| & \ge
|t|^{\frac{1}{j_0+1}}\cdot \frac{a}{2(r+1)CA^{r-1}}\nonumber \\
   &\ge
|t|^{\frac{1}{j_0+1}}\cdot \frac{d(q, q_j)}{2(r+1)CA^{r-1}}.
\end{align}
    Now assume  $|v|\le |t|^{\frac{1}{j_0+1}}\cdot
\frac{a}{(r+1)CA^{r-1}}$. In this case,
   \b{equation}\label{e7.25}
d(s, o_j)\ge \left|\tilde x_{j_0+1}+c_{j_0}(ad \,\tilde
x_1)^{j_0}(w'-(-1)^{j_0} tu)\right|^{\frac{1}{j_0+1}}.
\end{equation}
  Notice
$$w'-(-1)^{j_0}tu=w'+tu=2tu+v,$$
  so
   $$\tilde x_{j_0+1}+c_{j_0}(ad \,\tilde x_1)^{j_0}(w'-(-1)^{j_0} tu)=\tilde x_{j_0+1}+2tc_{j_0}(ad \,\tilde x_1)^{j_0}u +c_{j_0}(ad \,\tilde x_1)^{j_0}v.$$
Now   (\ref{e7.10})   and    (\ref{e7.20})  imply 
\b{align}\label{e7.50}
|c_{j_0}(ad \,\tilde x_1)^{j_0}v|  &\le |c_{j_0}|\cdot A^{j_0}\cdot |\tilde x_1|^{j_0}\cdot |v| \nonumber\\
 & \le C\cdot A^{j_0}\cdot \left(\frac{a}{101(r+1)CA^{r-1}}\right)^{j_0}\cdot |t|^{\frac{1}{j_0+1}}\cdot \frac{a}{(r+1)CA^{r-1}}\nonumber\\
  & \le  \frac{1}{101^{j_0}}\cdot \left(\frac{a}{r+1}\right)^{j_0+1}\cdot |t|^{\frac{1}{j_0+1}}.
\end{align}
  Since $|t|\ge 1$,    (\ref{e7.20}), (\ref{e7.22})  and  (\ref{e7.50})  imply
\b{align*} & |\tilde x_{j_0+1}+c_{j_0}(ad \,\tilde
x_1)^{j_0}(w'-(-1)^{j_0} tu)|\\
 & \ge\left(\frac{a}{r+1}\right)^{j_0+1}\cdot \left(|t|\cdot
(1-\frac{1}{(101CA^{r-1})^{j_0+1}})-\frac{1}{101^{j_0}}\cdot|t|^{\frac{1}{j_0+1}}-\frac{1}{(101CA^{r-1})^{j_0+1}}\right)\\
 &\ge
\frac{|t|}{2}\cdot\left(\frac{a}{r+1}\right)^{j_0+1}.
\end{align*}
Hence by (\ref{e7.25}) we have
    \b{equation}\label{e7.60}
d(s, o_j)\ge\frac{a}{(r+1) 2^{\frac{1}{j_0+1}}}\cdot
|t|^{\frac{1}{j_0+1}}\ge  \frac{d(q, q_j)}{(r+1)
2^{\frac{1}{j_0+1}}}\cdot |t|^{\frac{1}{j_0+1}}.
\end{equation}
  Now (\ref{e3.6.5}) follows from 
(\ref{e7.35})  and (\ref{e7.60}).  
The argument at the end of  the proof of Lemma 5.1 in \cite{X2}  now implies
 $$L_F(s)\le \frac{2(r+1)CA^{r-1} (\eta_1(1))^2}{|t|^{\frac{1}{j_0+1}}}\cdot
 L_F(q)\le \frac{2(r+1)CA^{r-1} (\eta_1(1))^2}{|t|^{\frac{1}{r}}}\cdot
 L_F(q).$$

\end{proof}

 The rest of the proof of Theorem \ref{abelianhigh}   is the same as in Section 5 of \cite{X2}.
\qed

The subspaces $W_1$ and $\hat W_1$ of $V_1$ are defined in Section \ref{invariant}. 

\b{Le}\label{l7.1} Let $\n=V_1\oplus \cdots\oplus V_r$ be a   
Carnot algebra.
\newline
  (1)   If  $r_1(\n)=1$,
then $[W_1, V_i]=0$ for all $i\ge 2$;\newline
 (2) If  $r_1(\n)\ge 2$  and      $r_{1, \C}(\n)\le 1$,  
then $[{\hat W}_1, V_i]=0$ for all $i\ge 2$.
\end{Le}

\b{proof} (1)  holds since $W_1$ is spanned by rank one elements and
$[X, V_i]=0$,    $i\ge 2$  for $X\in V_1$ with  $\text{rank}(X)=1$.

 (2)  Lemma \ref{crank1}
  implies 
$r_{1,\C}(\mathcal N)= 1$.  Hence  $W_{1, \C}\subset V_1\otimes \C$ is spanned by rank one
 elements.
   So  $[W_{1, \C}, \,V_i\otimes \C]=0$ for all $i\ge 2$.  It follows
   that   $[W_{1, \C},\, V_i]=0$ for all $i\ge 2$.
     By definition,  ${\hat W}_1$ is the set of real parts of
     elements in $W_{1, \C}$.  Hence (2) holds.

\end{proof}

Notice that Theorem \ref{main} is equivalent to the following:

\b{Th}\label{mainbody}
Let $\n$ be a non-rigid  Carnot algebra. If $\n$ is not one of the following  three classes of algebras, then
 it is quasisymmetrically rigid:\newline
 (1)     Euclidean algebras;\newline
 (2)  Heisenberg product    algebras;\newline
  (3)  complex  Heisenberg product     algebras.
\end{Th}

\b{proof}
By Theorem \ref{nonrex}  $\n$ has a reducible first layer.  If $\n$ is $2$-step, then 
Theorem \ref{2step1} implies $\n$ is
  quasisymmetrically rigid.   From now on we shall assume $\n$ is $r$-step with $r\ge 3$.     
  If $r_1(\n)=0$, then the claim follows from the main result in \cite{X3}.
   Assume  $r_1(\n)=1$. Then 
 Corollary  \ref{highrank}   implies 
 $W_1$   is  a  
non-trivial  proper subspace of $V_1$ invariant under  the action of
$\text{Aut}_g(\mathcal{N})$.    
If $W_1$ is abelian, then the claim follows from Lemma \ref{l7.1}(1)  and Theorem
  \ref{abelianhigh}.  If $W_1$ is not abelian, then the claim follows from Lemma \ref{r=1}.

Finally we assume $r_1(\n)\ge 2 $  and $r_{1, \C}(\n)\le 1$.   Then 
Corollary \ref{crank1cor}  
 implies 
 $\hat W_1$   is  a  
non-trivial  proper subspace of $V_1$ invariant under  the action of
$\text{Aut}_g(\mathcal{N})$.    
If $\hat W_1$ is abelian, then the claim follows from Lemma \ref{l7.1}(2)  and Theorem
  \ref{abelianhigh}.  If $W_1$ is not abelian, then the claim follows from Lemma \ref{r=2}.

\end{proof}

\section{Examples of biLipschitz   maps   on rank one groups}\label{examples}

In this section we construct an infinite dimensional space of
biLipschitz maps on Carnot groups whose Lie algebras  have  rank one
elements in the first layer.
  We remark that Ottazzi  \cite{O}
    has   previously  proved that  under the
  above assumption the space of smooth contact vector fields is
  infinite dimensional. The main point here is to construct explicit
  maps and also provide lots of non-smooth biLipschitz maps.

There are several different constructions that yield rank one Carnot algebras, see Section 4.2 of 
  \cite{X2}.

We shall need a result from \cite{X1}.
 The $n$-step ($n\ge 2$)  model Filiform algebra $\mathcal{F}^n$   is an $(n+1)$-dimensional  Lie
 algebra.
  It has a   basis $\{e_1, e_2, \cdots, e_{n+1}\}$
  and  the only non-trivial bracket relations are
    $[e_1, e_j]=e_{j+1}$  for $2\le j\le n$.
    The Lie algebra   $\mathcal{F}^n$   admits a  direct sum decomposition of vector subspaces
  $\mathcal{F}^n=V_1\oplus \cdots \oplus V_{n}$, where $V_1$ is the  linear subspace spanned by $e_1, e_2$, and $V_j$ ($2\le j\le n$)  is the linear subspace spanned by $e_{j+1}$.     It is easy to  check that
  $[V_1, V_j]=V_{j+1}$ for $1\le j\le n$,  where $V_{n+1}=\{0\}$.   Hence $\mathcal{F}^n$  is a
  graded      Lie algebra.
   The connected and simply connected Lie group with Lie algebra
$\mathcal{F}^n$    will be denoted by  $F^n$  and is called the
   $n$-step   model Filiform  group.

Let $h: \R\ra \R$ be a   $L$-Lipschitz   function  for some $L\ge
1$.
    Set $h_2=h$  and define $h_j:  \R  \ra \R$,  $3\le j\le n+1$, inductively as follows:
  $$h_j(x)=-\int_0^x h_{j-1}(s) ds.$$
    Define $G_h:   F^n\ra    F^n$   by:
    \b{equation}\label{F_h}
G_h\big(p)=p*\sum_{j=2}^{n+1}\big(h_j(x_1)\big)e_j
\end{equation}
    for   $p=\sum_{j=1}^{n+1}x_j e_j\in \mathcal F_n=F_n$. 
 It was  shown in  Section 3.4  of \cite{X1}   that  $G_h$ is $L'$-biLipschitz, where $L'$ is
 a constant depending only on $L$.

We also need the following Lemma of Ottazzi:
  \b{Le}\label{ottazzi} \e{(Lemma  3 in \cite{O})}
 Let $\n=V_1\oplus \cdots \oplus V_r$  be a Carnot algebra. Suppose
 that there exists $e_2\in V_1$   such that $\text{rank}(e_2)=1$.
 Then   $\n$ can be written as
  $$\n=\R  e_1\oplus \R e_2\oplus \R e_3\oplus  \cdots \oplus \R
  e_{n+1}\oplus {\mathcal  K},$$
    a vector space  direct sum, where  $e_1\in V_1$ and
     \b{equation}\label{e8.7}
[e_1, e_i]=e_{i+1}, \;\;2\le i\le n,
     \end{equation}
     \b{equation}
e_{n+1}\in { C}(\n),\;\; \; {\text{center of}} \;\; \n.
     \end{equation}
Moreover,
  $$\n_0=\R e_2\oplus \cdots \oplus \R e_{n+1}\oplus {\mathcal K}$$
   is an ideal   and
     $${\mathcal I}=\R e_2\oplus \cdots \oplus \R e_{n+1}$$
       is abelian.  Finally,     $V_1$ has a  basis $\{e_1,  e_2,
       U_1,  \cdots, U_s\}$ so that
\b{equation}\label{e8.10}
 [e_j,  U_i]=0,\; \forall j=2, \cdots, n+1, \forall i=1,
\cdots, s.
\end{equation}

  \end{Le}
\qed

Let 
 ${\mathcal F}=<e_1, e_2>=(\R   e_1\oplus \R  e_2)\oplus \R e_3 \oplus \cdots \oplus \R e_{n+1}$  be the subalgebra  of $\n$  generated by $e_1$ and $e_2$.
  Fix  a   vector space basis $\mathcal B$  of $\n$ such that each element of $\mathcal B$ belongs to some $V_k$
  and $e_i, U_j\in \mathcal B$ for all $i, j$.  
  For any $x\in \n$, we let $x_1$ be the coefficient of $e_1$ when $x$ is expressed as  a linear combination of elements in $\mathcal B$.  

\b{Le}\label{l8.20}
(1)   $[V_i, {\mathcal I}]=0$ for all $i\ge 2$;     \newline
  (2)     For any $x\in \n$  and any $y\in \mathcal I$, we have $[x, y]=x_1 [e_1, y]$;
  in particular,  $\mathcal I$ is an ideal of $\n$.
\end{Le}

\b{proof}
(1)     Notice that $W:=\R e_2$ satisfies the assumption of Lemma \ref{l7.11}.  
Now    the claim follows from  (\ref{e7.70})   for   $\tilde x_1=e_1$ and $u=e_2$.


(2) For any $x\in \n$  and any $y\in \mathcal I$,  the equality   $[x, y]=x_1 [e_1, y]$   follows from (1), (\ref{e8.10}) and the fact that $[e_2,  e_j]=0$ for all $j\ge 2$.    Then (\ref{e8.7}) implies 
  that   $\mathcal I$ is an ideal of $\n$.

\end{proof}

Here is the main result in this Section:

\b{Prop}\label{examplebilip}
Let $N$ be a Carnot group with Lie algebra
 $\n=V_1\oplus \cdots \oplus V_r$.
    Suppose there is some $e_2\in V_1$ with $\text{rank}(e_2)=1$.
Let  $U_1, \cdots, U_s\in V_1$  and $e_i\in V_{i-1}$, $2\le i\le n+1$ be given by Lemma \ref{ottazzi}.
  Let $\mathcal B$ be a basis of $\n$ such that each element of $\mathcal B$ belongs to some $V_k$
  and $e_i, U_j\in \mathcal B$ for all $i, j$.       
 Let $h:  \R \ra \R$ be a   Lipschitz function, and $h_j$ be as defined at the beginning of this section. 
  Then the map $F_h: N\ra N$,  
  $$F_h(x)=x*(\sum_{j=2}^{n+1}h_j(x_1) e_j)$$ 
   is biLipschitz, where $x_1$ is the   coefficient of  
    $e_1$  when  $x$  is written as a linear combination of elements of  $\mathcal B$.


\end{Prop}

\b{proof}    Let  each  $V_i$ be equipped with  the  inner product  such that $V_i\cap \mathcal B$ is an orthonormal basis.   Let   $d$ be  the corresponding homogeneous distance 
  on $\n=N$(see the last paragraph of Section \ref{basics}).

Let $x, \tilde x\in \n=N$ be  two arbitrary points.     We need to show that $d(x, \tilde x)$ and $d(F_h(x), F_h(\tilde x))$ are comparable.    Denote   $h_j=h_j(x_1)$  and $\tilde h_j=h_j(\tilde x_1)$. 
   Set  $y=\sum_{j=2}^{n+1} h_j e_j$ and $\tilde y=\sum_{j=2}^{n+1} \tilde h_j e_j$. 
  Then $F_h(x)=x*y$  and $F_h(\tilde x)=\tilde x*\tilde y$. 
 We have 
 $$d(x, \tilde x)=||(-x)*\tilde x||$$
     and 
 $$d(F_h(x),\, F_h(\tilde x))=d(x*y, \,  \tilde x*\tilde y)=||(-y)*(-x)*\tilde x*\tilde y||.$$

 By Lemma  \ref{l8.20}  (1),     $[V_i, {\mathcal I}]=0$  for all $i\ge  2$. 
  Notice $y, \tilde y\in \mathcal I$.  So Lemma \ref{cbhf} can be applied to obtain
 $$((-x)*\tilde x)*\tilde y=((-x)*\tilde x)  +\tilde y+\frac{1}{2}[(-x)*\tilde x,  \; \tilde y]+\sum_{j=2}^{r-1}c_j (ad\; (-x)*\tilde x)^j \tilde y.$$
  Since $((-x)*\tilde x)_1=\tilde x_1-x_1$,  
   Lemma   \ref{l8.20}   (2)  implies     
  \b{align}\label{e8.12}
((-x)*\tilde x)*\tilde y&=(-x)*\tilde x  +\tilde y+\frac{1}{2}(\tilde x_1-x_1)[e_1,   \tilde y]+\sum_{j=2}^{r-1}c_j(\tilde x_1-x_1)^j (ad\; e_1)^j \tilde y  \nonumber\\
&=(-x)*\tilde x  +\tilde y+\frac{1}{2}(\tilde x_1-x_1)[e_1,   \tilde y]+\sum_{j=2}^{n-1}c_j(\tilde x_1-x_1)^j (ad\; e_1)^j \tilde y.
\end{align}
 A similar argument applied to $(-y)*((-x)*\tilde x*\tilde y)$  yields
  \b{align*}
  &(-y)*((-x)*\tilde x*\tilde y)  \\
 &  =(-y)+   ((-x)*\tilde x*\tilde y)+
\frac{1}{2}[-y,  ((-x)*\tilde x*\tilde y)]+\sum_{j=2}^{r-1} (-1)^j c_j  (ad ((-x)*\tilde x*\tilde y))^j (-y)\\
  &=(-y)+ ((-x)*\tilde x*\tilde y)+
\frac{1}{2}[-y,  (\tilde x_1-x_1)e_1]+\sum_{j=2}^{n-1} (-1)^j  c_j (\tilde x_1-x_1)^j(ad\, e_1)^j (-y).
\end{align*}
  Using (\ref{e8.12}) we obtain:
  \b{equation}\label{e8.13}(- F_h(x ))* F_h(\tilde x)=
(-y)*((-x)*\tilde x*\tilde y)
  =(-x)*\tilde x +Y,
\end{equation}
      where  
 $$Y=(\tilde y-y)+\frac{1}{2}(\tilde x_1-x_1)[e_1,   \tilde y+y]
+\sum_{j=2}^{n-1}c_j (\tilde x_1-x_1)^j(ad\; e_1)^j (\tilde y-(-1)^jy).$$

Notice that   (\ref{e8.13}) implies   the coefficient of $e_1$  in both $(-x)*\tilde x$ and $(- F_h(x ))* F_h(\tilde x)$  are $\tilde x_1-x_1$.    Hence 
   \b{equation}\label{e8.14}
|\tilde x_1-x_1|\le ||(-x)*\tilde x||=d(x, \tilde x)
\end{equation}
 and 
  \b{equation}\label{e8.15}
|\tilde x_1-x_1|\le  d( F_h(x ),   F_h(\tilde x )).
\end{equation}

We observe that   the model Filiform group  $ F^n$  satisfies the assumption of the Proposition.   On  $\mathcal F^n=F^n$   we shall also  use the norm and homogeneous distance 
 defined at the beginning of the proof. All the above calculations and estimates are valid for  $\mathcal F^n$ as well.  
Notice that $ F_h(x_1e_1)=G_h(x_1e_1)$  and $ F_h(\tilde x_1e_1)=G_h(\tilde x_1e_1)$.
  Since $(-x_1e_1)*\tilde x_1 e_1=(\tilde x_1-x_1) e_1$,
   (\ref{e8.13}) implies 
   $$(-F_h(x_1 e_1))*F_h(\tilde x_1 e_1)=(\tilde x_1-x_1) e_1  +Y
=(- G_h(x_1 e_1))* G_h(\tilde x_1 e_1).$$
  Since  $G_h$ is   $L$-biLipschitz for some $L\ge 1$,  the triangle inequality now implies
   \b{align}\label{e8.16}
||Y||\le ||(\tilde x_1-x_1) e_1  +Y||&=d(G_h(x_1 e_1),  G_h(\tilde x_1 e_1))\le L\cdot d(x_1 e_1, \tilde x_1 e_1)\nonumber \\
 &  =L\cdot |(-x_1 e_1)* \tilde x_1 e_1|=L\cdot |\tilde x_1-x_1|.
\end{align}
  Now (\ref{e8.13}),  triangle inequality,    (\ref{e8.16})  and   (\ref{e8.14})  imply 
 \b{align*}
d( F_h(x ),   F_h(\tilde x )) & =||(- F_h(x ))* F_h(\tilde x )||\\
  & =||(-x)*\tilde x +Y||\le 
||(-x)*\tilde x||+||Y||\le d(x, \tilde x)+L\cdot |\tilde x_1-x_1|\\
  & \le  (L+1)\cdot d(x, \tilde x).
\end{align*}
Similarly,  using (\ref{e8.15}) instead of  (\ref{e8.14})  we obtain 
\b{align*}
d(x, \tilde x)=||(-x)*\tilde x||= &  ||(- F_h(x ))* F_h(\tilde x )-Y|| \\
& \le 
||(- F_h(x ))* F_h(\tilde x )|| +||Y||\\
 &  \le  d( F_h(x ),   F_h(\tilde x ))   
+  L\cdot |\tilde x_1-x_1|\\
 & \le d( F_h(x ),   F_h(\tilde x ))   
+  L\cdot d( F_h(x ),   F_h(\tilde x ))  \\
 & =(L+1)d( F_h(x ),   F_h(\tilde x ))   .
\end{align*}
 Hence $ F_h$ is biLipschitz.


\end{proof}

 \addcontentsline{toc}{subsection}{References}

\noindent Address:

\noindent Xiangdong Xie: Dept. of Mathematics  and   Statistics,   Bowling Green  State  University,
  Bowling Green,  OH,   U.S.A.\hskip .4cm E-mail:   xiex@bgsu.edu

\end{document}